\documentclass{amsart}

\usepackage[a4paper, margin = 60pt]{geometry}
\usepackage{amsmath,amsthm,amssymb,bm,mathrsfs}
\usepackage{longtable,array}
\usepackage{xr-hyper}
\usepackage{hyperref}
\usepackage{cleveref}
\usepackage{graphicx,color}
\usepackage{tikz}
\usepackage[T1]{fontenc}
\IfFileExists{tgtermes.sty}{\usepackage{tgtermes}}{}
\numberwithin{equation}{section}

\newtheorem{thm}{Theorem}[section]
\newtheorem{lem}[thm]{Lemma}
\newtheorem{prop}[thm]{Proposition}
\newtheorem{cor}[thm]{Corollary}
\theoremstyle{definition}
\newtheorem{exam}[thm]{Example}
\newtheorem{defn}[thm]{Definition}

\newtheorem{remark}[thm]{Remark}

\crefname{lem}{Lemma}{Lemmas}
\crefname{thm}{Theorem}{Theorems}
\crefname{prop}{Proposition}{Propositions}
\crefname{question}{Question}{Questions}
\crefname{defn}{Definition}{Definitions}
\crefname{conj}{Conjecture}{Conjectures}
\crefname{figure}{Figure}{Figures}
\crefname{cor}{Corollary}{Corollaries} 
\crefformat{equation}{(#2#1#3)}
\Crefformat{equation}{Equation #2(#1)#3}

\AddToHook{env/lem/begin}{\crefalias{thm}{lem}}
\AddToHook{env/prop/begin}{\crefalias{thm}{prop}}
\AddToHook{env/cor/begin}{\crefalias{thm}{cor}}
\AddToHook{env/exam/begin}{\crefalias{thm}{exam}}
\AddToHook{env/defn/begin}{\crefalias{thm}{defn}}
\AddToHook{env/conj/begin}{\crefalias{thm}{conj}}
\AddToHook{env/question/begin}{\crefalias{thm}{question}}
\AddToHook{env/problem/begin}{\crefalias{thm}{problem}}
\AddToHook{env/remark/begin}{\crefalias{thm}{remark}}

\newcommand\lc{\operatorname{lc}}
\newcommand\wf{\widetilde{f}}

\newcommand\A{\mathcal{A}}
\newcommand\D{\mathcal{D}}
\newcommand\AW{\mathscr{D}_q}
\newcommand\Rq{\mathcal{R}_q}
\newcommand\ii{\mathrm{i}}
\newcommand\Deltaop{\Delta}

\newcommand\res{\operatorname{res}}
\newcommand\disc{\operatorname{disc}}
\newcommand\OPS{\operatorname{OPS}}

\newcommand\LL{\mathcal{L}}

\renewcommand\vec[1]{{\bm #1}}

\definecolor{darkblue}{rgb}{0.0,0,0.7}
\renewcommand\emph[1]{\textcolor{darkblue}{\it #1}}

\title{Discriminants and resultants of orthogonal polynomials}
\author{Mourad E. H. Ismail}
\address{Department of Mathematics, University of Louisiana at Lafayette, Lafayette, LA 70504-3568, USA}
\email{mourad.eh.ismail@gmail.com}

\author{Jihyeug Jang}
\address{Department of Mathematics and Institute of Mathematical Sciences, College of Natural Sciences, Pusan National University, Busan, Republic of Korea}
\email{jihyeugjang@pusan.ac.kr}

\author{Jang Soo Kim}
\address{Department of Mathematics,
Sungkyunkwan University (SKKU), Suwon, Gyeonggi-do 16419, South Korea}
\email{jangsookim@skku.edu}

\author{Dennis Stanton}
\address{School of Mathematics,
University of Minnesota, Minneapolis, Minnesota 55455, USA}
\email{stanton@math.umn.edu}
\begin{document}

\begin{abstract}
The concept of a linked pair of two sequences of polynomials, one of which is orthogonal, is introduced. Such a pair always has a resultant that factors. Two methods for producing such linked pairs are given: one for general orthogonal polynomials and the other for any polynomial in the Askey scheme. Generalized discriminants are introduced and evaluated, including a general product formula for each operator in the Askey scheme. A new notion of companion orthogonal polynomials is also introduced and used to evaluate the discriminants of orthogonal polynomials. These methods recover all known formulas for the generalized discriminants of classical polynomials and provide new ones.
\end{abstract}

\maketitle
\tableofcontents

\section{Introduction}
\label{sec:introduction}

In the 1880's Hilbert \cite{Hilbert1888} and Stieltjes
\cite{Stieltjes1885b,Stieltjes1885a} found the discriminant of Jacobi
polynomials. The discriminant of a polynomial is a multiple of the
resultant of two polynomials: the polynomial and its derivative. Thus
Hilbert and Stieltjes had found a result for the resultant of \(p\)
and \(Dp\), denoted \( \res(p,Dp) \). Analogous results follow for
Laguerre and Hermite polynomials. Resultants are used in algebraic
geometry, elimination theory, and number theory to detect common roots
and eliminate variables.

Stieltjes was studying an electrostatics problem whose equilibrium
positions have particles located at the zeros of Jacobi polynomials.
The energy to be minimized was related to the discriminant of these
polynomials. This problem was later generalized in many directions by
several researchers. For example, the electrostatics problem can be
found in \cite[Section~3.5]{Ismail} and \cite{Szego1975}. Other electrostatics
versions are given by Ismail; see \cite{Ismail2000,Ismail2000a}.

Discriminants of orthogonal polynomials have been studied by many
researchers; see for example
\cite{Dilcher2005,SawaUchida2019,IsmailWang2019,Matsumura2024}. Ismail
generalized the discriminant evaluation in two ways. First, by
considering orthogonal polynomials whose measure has certain
properties, he gave a general formula for the discriminant
\cite{Ismail1998}; see also \cite[Theorem~3.4.2]{Ismail}. Second, he
generalized the discriminant to \( \res(p,Tp) \) by replacing the
derivative operator \( D \) by a degree-lowering operator \( T \) on
polynomials; see \cite[(6.4.2)]{Ismail}. Considering different choices
of \( T \) naturally leads to discrete and \(q\)-versions of the
discriminant; see \cite{Ismail2003,Ismail2004} and \cite[(6.4.3),
(18.6.11)]{Ismail}. Some of the polynomials in the Askey and
\(q\)-Askey schemes then have explicit product formulas for these
discriminants \cite[Sec.~6.4, Theorem~15.5.3]{Ismail}. However, not
all polynomials in the Askey and \(q\)-Askey schemes had generalized
discriminants that factored; for example, the Askey--Wilson
polynomials did not. This paper extends these results.

The key idea is a mixed
three-term relation that connects two polynomial sequences, one of which
is orthogonal; we call such a pair a \emph{linked pair}. There is always a resultant that
factors for linked pairs. This motivated us to further generalize
discriminants as \( \res(\A_{u,v} p, \D_{u,v} p) \) using divided
difference operators \( \D_{u,v} \) and averaging operators
\( \A_{u,v} \), where \( u \) and \( v \) are functions satisfying
certain conditions. In the classical discriminant case, the divided
difference operator is the derivative and the averaging operator is
the identity operator. If \( v \) is the identity function, then our
discriminant \( \res(\A_{u,v} p, \D_{u,v} p) \) reduces to Ismail's
discriminant \( \res(p, \D_{u,v} p) \) up to scalar multiplication.
When the divided difference operator is nice, we give an explicit
product formula for this new generalized discriminant of a polynomial
in terms of its zeros.

We give two methods for finding linked pairs. First, suppose that
\( \{p_n(x)\} \) is a monic orthogonal polynomial sequence (OPS) such
that \( \{\D_{u,v}p_n\} \) is also an OPS. We show that up to scalar
multiplication \( \{\A_{u,v}p_n\} \) and \( \{\D_{u,v}p_n\} \) form a
linked pair. In this case, we provide an explicit product formula for
the generalized discriminant \( \res(\A_{u,v} p_n, \D_{u,v} p_n) \)
using their three-term recurrence coefficients. This method applies to
every polynomial family in the Askey and \(q\)-Askey schemes, hence it
has a corresponding resultant and the generalized discriminant that
factor.

The second method applies to general orthogonal polynomials: given
orthogonal polynomials, we can construct a corresponding polynomial
sequence to form a linked pair. This is accomplished by multiplying
the measure by a polynomial of degree at most two. The classical
discriminant for Jacobi polynomials follows from this method because
the derivative of a Jacobi polynomial is orthogonal with respect to
the measure obtained by multiplying the Jacobi measure by the
quadratic polynomial \((1-x)(1+x)\). This general result can be
applied to the classical polynomial families, including the
Askey--Wilson polynomials. All known generalized discriminants for the
orthogonal polynomials in the Askey and \(q\)-Askey schemes can be
proved using at least one of these two methods.

We also introduce a new class of orthogonal polynomial sequences,
called \emph{\( (f,g) \)-companion orthogonal polynomial sequences},
or \( (f,g) \)-COPS, associated with polynomials \( f \) and \( g \).
In a special case, the moments of the \( (f,f') \)-COPS associated
with a polynomial \(f\) are equal to the power-sum symmetric
polynomials in the zeros of \(f\), divided by \(\deg f\). This allows
ordinary discriminants to factor naturally because the discriminant of
a monic polynomial is the Hankel determinant of the power-sum
symmetric polynomials of its zeros and the Hankel determinant of the
moments of an OPS factors into its recurrence coefficients. This
observation was our original starting point for investigating the
factorization phenomena exhibited by the discriminants of various
orthogonal polynomials.

The paper is organized as follows. Basic definitions and a review of
resultants, discriminants, and the Hahn difference operator are given
in \Cref{sec:prelim}.
In \Cref{sec:hahn-diff-oper}, we introduce the generalized
discriminant with two operators. If the two operators are a
divided difference operator and its associated averaging operator, a
product formula for the general discriminant in terms of the zeros of
a polynomial is given in
\Cref{pro:disc-prod2}. This generalizes the classical discriminant of a
monic polynomial, which is the square of Vandermonde's determinant of the zeros. The
Askey--Wilson operator is a special case of this divided difference
operator.
Section~\ref{sec:mixed-three-term} defines linked pairs, and then they
are used to prove the factorization formulas; for resultants see
\Cref{thm:linked-pair-resultant} and for generalized discriminants see
Corollary~\ref{cor:avg-T-disc}. If the divided difference operator
maps an OPS to another OPS, then we explicitly obtain a linked pair in
\Cref{thm:mixed-3-term}. We also show that \Cref{thm:mixed-3-term}
applies to all orthogonal polynomials in the Askey and \(q\)-Askey
schemes.
In \Cref{sec:result-orth-polyn}, we prove the general results to find
linked pairs, \Cref{linearthm2} and \Cref{quadthm}, by multiplying
the measure by \(x-A\) or \((x-A)(x-B)\).
In \Cref{sec:companion-ops}, we introduce COPS and find an algorithm
to construct them.

The results of this paper are based on Schur's fundamental result,
Lemma~\ref{thm:schur}.

% \subsection*{Acknowledgments}

\begingroup
\renewcommand{\addcontentsline}[3]{}% Omit this starred subsection from the amsart TOC.
\subsection*{Declaration on Generative AI}

The ChatGPT 5.5, 5.6, and 6 Pro models were used for the following
tasks:
\begin{itemize}
\item Finding and proving \Cref{pro:disc-prod2} and \Cref{prop:res-hankel}.
\item Searching references, proofreading the manuscript, and
  conducting computational experiments with Sage.
\end{itemize}
All proofs obtained with the assistance of these AI tools have been
independently verified and rewritten by the authors.
\par\endgroup

\section{Preliminaries}
\label{sec:prelim}

This section contains definitions and the relevant notation used
throughout the paper.

\subsection{Basic definitions}

We begin with the definition of orthogonal polynomials. We use the
convention that allows us to treat both infinite and finite orthogonal
polynomial sequences at the same time.

\begin{defn}
  Let \( \LL \) be a linear functional on the space of polynomials and
  let \( N \in \{1,2,\dots\}\cup\{\infty\} \). A sequence
  \( \{f_n(x):0\le n<N \} \) of polynomials with \(\deg f_n=n\) is
  called an \emph{orthogonal polynomial sequence (OPS)} with respect
  to \( \LL \) if \( \LL(f_n(x)f_m(x))=0 \) for \( 0\le n\ne m <N \),
  and \( \LL(f_n(x)^2)\ne 0 \) for \( 0\le n <N \). It is called a
  \emph{monic OPS} if every \(f_n\) is monic.
\end{defn}

Favard's theorem states that \( \{f_n(x):0\le n<N \} \) is a monic OPS
with respect to some linear functional \( \LL \) if and only if, for
some coefficients \(b_n\) and \(\lambda_n\), it satisfies
  \[
    f_{n+1}(x)=(x-b_n)f_n(x)-\lambda_nf_{n-1}(x)
    \qquad (0\le n<N-1),
  \]
  where \(f_{-1}(x)=0\), \(f_0(x)=1\), and \(\lambda_n\ne0\) for
  \(1\le n<N-1\). In this case, we write
\begin{equation}\label{eq:OPS}
  \{f_n(x)\} = \OPS(\{b_n\},\{\lambda_n\}).
\end{equation}

In this paper, we consider monic orthogonal polynomials in the Askey
and \(q\)-Askey schemes. We give only the definition of the monic
Askey--Wilson polynomials; see \cite{KLS} for the definitions of the
other families. For \(q\)-series and the basic hypergeometric series
\( {}_r\phi_s \), we follow the notation of \cite{gasper}.

The monic \emph{Askey--Wilson polynomial} is defined by
\[
  p_n(x;a,b,c,d|q)
  =
  \frac{a^{-n}(ab,ac,ad;q)_n}
       {2^n(abcdq^{n-1};q)_n}
  {}_4\phi_3\!\left(
    \begin{matrix}
      q^{-n},abcdq^{n-1},az,az^{-1}\\
      ab,ac,ad
    \end{matrix};q,q
  \right),
\]
where \(x=(z+z^{-1})/2\). They satisfy the three-term recurrence
relation
\[
  p_{n+1}(x;a,b,c,d|q)=(x-b_n)p_n(x;a,b,c,d|q)-\lambda_n p_{n-1}(x;a,b,c,d|q),
\]
where \( p_0(x;a,b,c,d|q)=1 \), \( p_1(x;a,b,c,d|q)=x-b_0 \), and
\begin{align}
\label{eq:AW-b,la}  b_n&=\frac12\bigl(a+a^{-1}-A_n-C_n\bigr),
  \qquad \lambda_n=\frac14A_{n-1}C_n,\\
  \notag A_n&=
  \frac{(1-abq^n)(1-acq^n)(1-adq^n)(1-abcdq^{n-1})}
       {a(1-abcdq^{2n-1})(1-abcdq^{2n})},\\
  \notag C_n&=
  \frac{a(1-q^n)(1-bcq^{n-1})(1-bdq^{n-1})(1-cdq^{n-1})}
       {(1-abcdq^{2n-2})(1-abcdq^{2n-1})}.
\end{align}

An \emph{operator} is a linear transformation on the space of
polynomials. In this paper, we consider various operators. Three basic
operators are the \emph{derivative operator} \( D \),
the \emph{finite difference operator} \( \Delta \), and the
\emph{\( q \)-difference operator} \( D_q \), which are defined by
\( D=d/dx \) and 
\begin{align*}
  \Delta f(x) &= f(x+1) - f(x),\\
    D_q f(x) &= \frac{f(x)-f(qx)}{(1-q)x}.
\end{align*}

\subsection{Resultants and discriminants}

We recall the definitions of resultants and discriminants and their
basic properties.

\begin{defn}\label{def:resultant}
 Given two polynomials \(f(x),g(x)\) of degree \(n\) and \(m\) with
\[
f(x)=A\prod_{j=1}^n (x-\alpha_{j}), \quad g(x)=B\prod_{j=1}^m (x-\beta_{j}),
\]
the \emph{resultant} \( \res(f(x),g(x)) \) of \( f(x) \) and
\( g(x) \) is defined by
\[
 \res(f(x),g(x))= A^mB^n \prod_{j=1}^n\prod_{k=1}^m (\alpha_j-\beta_k)
= A^m\prod_{j=1}^n g(\alpha_j)=(-1)^{mn}B^n  \prod_{k=1}^m f(\beta_k).
\]
\end{defn}

We will use the following standard properties of resultants, which
follow easily from the definition.

\begin{prop}\label{pro:4}
  Let \( f \), \( g \), and \( h \) be polynomials.
  \begin{enumerate}
  \item We have \( \res(f,g)=0 \) if and only if \( f \) and \( g \)
    have a common zero.
  \item The resultant is multiplicative in each argument:
\[
  \res(f,gh)=\res(f,g)\res(f,h).
\]
\item 
If \(\deg(f+hg) = \deg f\), then
\[
  \res(f+hg,g)=\res(f,g).
\]
  \end{enumerate}
\end{prop}

We now define discriminants, which are closely related to resultants.

\begin{defn}\label{def:6}
The \emph{discriminant} \( \disc(f) \) of the polynomial
\( f(x) = A(x-\alpha_1)\cdots (x-\alpha_n) \) is defined by
\[
  \disc(f)
  = A^{2n-2} \prod_{1\le i<j\le n}\left(\alpha_i-\alpha_j\right)^2
  =(-1)^{\binom{n}{2}} A^{2n-2} \prod_{1\le i\ne j\le n}\left(\alpha_i-\alpha_j\right).
\]
\end{defn}

Since
\[
  \res(f,Df) = A^{n-1} \prod_{j=1}^{n} Df(\alpha_j)
  = A^{2n-1} \prod_{1\le i\ne j\le n} (\alpha_i-\alpha_j),
\]
we have
\[
\disc(f)= (-1)^{\binom{n}{2}} A^{-1} \res(f,Df).
\]
Ismail \cite{Ismail2000} introduced the concept of a discriminant with respect to an operator.

\begin{defn}\label{def:disc}
  Let \( f(x) = A(x-\alpha_1)\cdots (x-\alpha_n) \) be a
  polynomial of degree \( n \), and suppose that \( Tf \) has degree
  \( m \). The \emph{discriminant}
  \( \disc(f;T) \) of \( f \) with respect to an operator \(T\) is
  defined by
\[
\disc(f;T)=(-1)^{\binom{n}{2}} A^{-1} \res(f,Tf)=(-1)^{\binom{n}{2}} A^{m-1} \prod_{j=1}^n(T f)\left(\alpha_j\right).
\]
\end{defn}

The following Schur's lemma shows that the resultant of two
consecutive monic orthogonal polynomials has a simple product formula;
see \cite[p. 67, Lemma 3.4.1]{Ismail}. This is a main tool to prove
various product formulas in this paper.

\begin{lem}[Schur's lemma]
  \label{thm:schur} 
  Let \( \{p_n(x)\} = \OPS(\{b_n\},\{\lambda_n\}) \). Then 
$$
\res(p_{n+1}(x), p_n(x))= \res(p_{n}(x), p_{n+1}(x))= (-1)^ {\binom{n+1}{2}} 
\prod_{j=1}^n \lambda_{j}^{j}.
$$
\end{lem}

\subsection{Hahn difference operator}

The Hahn difference operator \cite{Hahn1949} is a useful tool that
contains many important operators; see also
\cite[Definition~2.5]{Hamza2015}. To define this, we use the following
convention for writing a quotient. For a polynomial \( f=f(x) \) and
any functions \( u=u(x) \) and \( v=v(x) \), we write
\begin{equation}\label{eq:u-v}
  \frac{f(u)-f(v)}{u-v} = \frac{f(y)-f(z)}{y-z}\bigg|_{y=u, z=v}.
\end{equation}
Here the right-hand side means we substitute \( y=u \) and \( z=v \)
into \( (f(y)-f(z))/(y-z) \), which is a polynomial in the
indeterminates \( y \) and \( z \). Note that this is well defined for
any functions \( u \) and \( v \), including the case \( u=v \).

\begin{defn} \label{def:1} Let \( u=u(x)=ax+b \) be
  an affine polynomial in \( x \). The \emph{Hahn difference operator}
  \( \D_u \) with respect to \( u \) is the operator on the
  space of polynomials defined by
\[
  \D_u(f(x)) = \frac{f(u(x)) - f(x)}{u(x) - x}.
\]
\end{defn}

Note that the condition that \( u \) is a polynomials in \( x \) is to
ensure that \( \D_u \) is an operator on polynomials. By the convention
\eqref{eq:u-v}, if \( u(x)=x \), then
\[
  \D_u(x^n) = \left.\frac{y^n-z^n}{y-z} \right|_{y=x,z=x}
   = \left.\left( \sum_{i=0}^{n-1} y^iz^{n-1-i}  \right) \right|_{y=x,z=x}
   = nx^{n-1}.
\]
Hence, in this case, \( \D_u \) is the same as the derivative operator
\( D \).

The discriminant \( \disc(f(x);\D_u) \) has a
product formula in terms of the zeros of \( f(x) \); see
\cite[(6.4.3), (18.6.11)]{Ismail}.

\begin{prop}\label{pro:disc-prod}
  Let \( u \) be an affine polynomial in \( x \). For a monic polynomial
  \( f(x)=\prod_{k=1}^n(x-\alpha_k) \), we have
  \[
    \disc(f;\D_u)
    =\prod_{1\le i<j\le n}
      (\alpha_i-u(\alpha_j))(u(\alpha_i)-\alpha_j).
  \]
\end{prop}

\begin{proof}
  Since \( f(\alpha_i)=0 \), we have
  \[
    (\D_u f)(\alpha_i)
    =\frac{f(u(\alpha_i))-f(\alpha_i)}{u(\alpha_i)-\alpha_i}
    =\prod_{j\ne i}(u(\alpha_i)-\alpha_j).
  \]
  Therefore, by \Cref{def:disc}, we obtain
  \[
    \disc(f;\D_u)
    =(-1)^{\binom n2}
      \prod_{i=1}^n\prod_{j\ne i}(u(\alpha_i)-\alpha_j)
          =\prod_{1\le i<j\le n}
      (\alpha_i-u(\alpha_j))(u(\alpha_i)-\alpha_j). \qedhere
  \]
\end{proof}

\begin{exam}\label{exa:2}
  Let \( f(x) = (x-\alpha_1)\cdots (x-\alpha_n) \).
  \begin{enumerate}
  \item If \( u(x) = x \), then \( \D_u = D \). In this
    case, \( \disc(f;\D_u) \) is the usual discriminant:
    \[
      \disc(f;\D_u) = \prod_{1\le i<j\le n} (\alpha_i-\alpha_j)^2.
    \]
  \item If \( u(x) = x+1 \), then \( \D_u \) is the difference
    operator \( \Delta \). In this case, we have
    \[
      \disc(f;\D_u) = \prod_{1 \leq i<j \leq n}
      \left(\alpha_i-\alpha_j-1\right)\left(\alpha_i-\alpha_j+1\right) .
    \]
  \item If \( u(x) = qx \), then \( \D_u \) is the
    \( q \)-difference operator \( D_q \). In this case, we have
    \[
      \disc(f;\D_u) = \prod_{1 \leq i<j \leq n}
      \left(\alpha_i-q\alpha_j\right)\left(q\alpha_i-\alpha_j\right) .
    \]
  \end{enumerate}
\end{exam}

\section{Divided difference operators and discriminants}
\label{sec:hahn-diff-oper}

In this section, we introduce the divided difference
operator \( \D_{u,v} \) and the averaging operator \( \A_{u,v} \). We
review the forward shift operators of the orthogonal polynomials in
the Askey and \( q \)-Askey schemes, and describe them using the
divided difference operator in a uniform way. We also
introduce the generalized discriminant \( \disc(f;\A_{u,v},\D_{u,v}) \)
associated with these two operators. We show that it has a product
formula in terms of the zeros of the polynomial \( f \).

\subsection{Divided difference operator}

We begin with the definition. Recall the convention \eqref{eq:u-v}.

\begin{defn}\label{def:Hahn2}
  Let \( u \) and \( v \) be functions in \( x \) such that \( u+v \)
  and \( uv \) are polynomials in \( x \). For a polynomial \( f \),
  define the \emph{divided difference operator}
  \( \D_{u,v} \) and the \emph{averaging operator} \( \A_{u,v} \) by
\[
  \D_{u,v}f
  = \frac{f(u) - f(v)}{u-v},
  \qquad 
  \mathcal{A}_{u,v}f
  = \frac{f(u) + f(v)}2.
\]
If \( T=\D_{u,v} \), we denote by
\( \A_{T} = \A_{u,v} \) the associated averaging operator.
\end{defn}

Note that if \( u(x)=ax+b \) and \( v(x)=x \), then \( \D_{u,v} \) reduces to the Hahn
operator \( \D_u \). By the convention \eqref{eq:u-v}, if \( u=v \),
then \( \D_{u,u}f=(Df)\circ u \). In particular, if \( u(x)=x \), then
\( \D_{u,u} \) becomes the derivative operator \( D \).

The condition that \( u+v \) and \( uv \) are
polynomials in \( x \) is equivalent to the condition that the
following are polynomials in \( x \) for all integers \( n\ge0 \):
\[
  \A_{u,v}(x^n)=\frac{u^n+v^n}{2},
  \qquad
  \D_{u,v}(x^n)=\frac{u^n-v^n}{u-v}.
\]
Indeed, since the above quantities are symmetric polynomials in
\( u \) and \( v \), they can be computed from \(U=u+v\) and \(V=uv\).
To see this, let
\[
  S_m=u^m+v^m, \qquad H_m=\frac{u^m-v^m}{u-v}.
\]
Then \( H_m \) and \( S_m \) can be computed using the recurrences:
\begin{equation}\label{eq:SH-rec}
  S_m=US_{m-1}-VS_{m-2}, \qquad
  H_m=UH_{m-1}-VH_{m-2}\qquad (m\ge2),
\end{equation}
where \( S_0=2 \), \( S_1=U \), \( H_0=0 \), and \( H_1=1 \).
Thus, if \(f(x)=\sum_{m=0}^n c_mx^m\), then
\begin{equation}\label{eq:HS}
  \A_{u,v}f=\frac12\sum_{m=0}^n c_mS_m,\qquad
  \D_{u,v}f=\sum_{m=0}^n c_mH_m.
\end{equation}

The operators \( \D_{u,v} \) and \( \A_{u,v} \) have the following
product rules; see \cite[(12.1.22)]{Ismail} and
\cite[(2.5)]{MbounaSuzuki2022} for the symmetric product rules
\eqref{eq:avg-product} and \eqref{eq:A-product}, respectively, in the
case of the Askey--Wilson operator.

\begin{prop}\label{lem:product-rule}
  Let \( u \) and \( v \) be functions in \( x \) such that \( u+v \)
  and \( uv \) are polynomials in \( x \). For all polynomials \( f \)
  and \( g \), we have
  \begin{align}\label{eq:leibniz}
    \D_{u,v}(fg)
    &= (f\circ u)\,(\D_{u,v}g) + (\D_{u,v}f)(g\circ v),\\
    \label{eq:avg-product}
    \D_{u,v}(fg)
    &= (\mathcal{A}_{u,v}f)\,(\D_{u,v}g)
      + (\D_{u,v}f)\,(\mathcal{A}_{u,v}g),\\
    \label{eq:A-product}
    \A_{u,v}(fg)
    &=
    (\A_{u,v}f)(\A_{u,v}g)
    +\frac{(u-v)^2}{4}(\D_{u,v}f)(\D_{u,v}g).
  \end{align}
\end{prop}

\begin{proof}
  If \( u=v \), then \( \D_{u,v}f=(Df)\circ u \) and
  \( \A_{u,v}f=f\circ u \). Hence, the first two identities follow
  from the ordinary product rule, evaluated at \(u(x)\), and the third
  identity is immediate. Suppose now that
  \( u\ne v \). Let \(a=f\circ u\), \(b=f\circ v\),
  \(c=g\circ u\), and \(d=g\circ v\). Then
  \(\D_{u,v}(fg)=\frac{ac-bd}{u-v}\),
  \(\A_{u,v}(fg)=\frac{ac+bd}{2}\),
  \(\D_{u,v}f=\frac{a-b}{u-v}\), 
  \(\D_{u,v}g=\frac{c-d}{u-v}\),
  \(\A_{u,v}f=\frac{a+b}{2}\), and
  \(\A_{u,v}g=\frac{c+d}{2}\).
  The first identity follows from
  \[
    \D_{u,v}(fg)
   = \frac{ac-bd}{u-v}=\frac{a(c-d)+(a-b)d}{u-v}
     = (f\circ u)\,(\D_{u,v}g) + (\D_{u,v}f)(g\circ v).
  \]
  Similarly, we obtain the second identity from
  \[
    \D_{u,v}(fg) =
    \frac{ac-bd}{u-v}=\frac{\frac{a+b}{2}(c-d)+(a-b)\frac{c+d}{2}}{u-v}
    = (\mathcal{A}_{u,v}f)\,(\D_{u,v}g) + (\D_{u,v}f)\,(\mathcal{A}_{u,v}g).
  \]
  Finally,
  \[
    (\A_{u,v}f)(\A_{u,v}g)
    +\frac{(u-v)^2}{4}(\D_{u,v}f)(\D_{u,v}g)
    =
    \frac{(a+b)(c+d)+(a-b)(c-d)}4  
    =
    \frac{ac+bd}{2}
    =
    \A_{u,v}(fg). \qedhere
  \]
\end{proof}

For our purposes, we impose one additional condition beyond the
requirement that \( u+v \) and \( uv \) are polynomials.

\begin{defn}\label{def:2}
Let \( u \) and \( v \) be functions of \( x \).
We say that \( (u,v) \) is a \emph{degree-lowering pair}
if the following conditions hold:
\begin{itemize}
\item Both \( u+v \) and \( uv \) are polynomials in \( x \).
\item If \( f \) is a polynomial of degree \( n\ge1 \), then
  \( \D_{u,v}f \) is a polynomial of degree \( n-1 \).
\end{itemize}
\end{defn}

Suppose that \( (u,v) \) is a degree-lowering pair. By
\eqref{eq:avg-product} with \( f=g=x^n \), we have
\( \D_{u,v}x^{2n} =2(\A_{u,v}x^n)(\D_{u,v}x^n) \). Since \(\D_{u,v}\)
lowers degree by one, we have
\[
  \deg \A_{u,v}x^n=n \qquad (n\ge1).
\]
In particular, \( \A_{u,v}x=\rho x+\kappa \) for some constants
\( \rho \) and \( \kappa \) with \( \rho\ne0 \).

From now on, we will only consider operators \( \D_{u,v} \) for
degree-lowering pairs \( (u,v) \). Such functions \( u \) and \( v \)
can be written explicitly as follows.

\begin{prop}\label{pro:good-pair}
  Let \( (u,v) \) be a degree-lowering pair. Then there exist unique polynomials
  \( U=U(x) \) and \( V=V(x) \) with \( \deg U=1 \) and
  \( \deg V\le 2 \) such that \( u \) and \( v \) are the roots of
  \( t^2-Ut+V=0 \), that is,
  \[
    \{u,v\}
    =
    \left\{
      \frac{U+\sqrt{U^2-4V}}{2},
      \frac{U-\sqrt{U^2-4V}}{2}
    \right\}.
  \]
\end{prop}

\begin{proof}
  Let \( U=u+v \) and \( V=uv \). Then \( u \) and \( v \) are the
  roots of \( t^2-Ut+V=0 \). Since \( (u,v) \) is a degree-lowering pair, \( U \)
  and \( V \) are polynomials. Since
  \begin{align*}
    \D_{u,v}x^2&=u+v=U,\\
    \D_{u,v}x^3
    &=\frac{u^3-v^3}{u-v} =u^2+uv+v^2 =U^2-V,
  \end{align*}
  we have \( \deg U=1 \) and \( \deg(U^2-V)=2 \), which implies
  \( \deg V\le 2 \). The uniqueness of \( U \)
  and \( V \) follows from the above equations.
\end{proof}

\subsection{Forward shift operators in Askey schemes}

Each OPS \( \{p_n(x;\vec t)\} \) with parameters \( \vec t \) in the
Askey and \( q \)-Askey schemes has a \emph{forward shift operator},
which is an operator \( T \) such that
\( T p_n(x;\vec t) = c_n p_{n-1}(r x+s;\vec t') \), where \( c_n \),
\(r\ne0\) and \(s\) are constants, and \( \vec t' \) denotes the
parameters shifted forward. For example, the Askey--Wilson operator
\( \AW \) is the forward shift operator for the monic Askey--Wilson
polynomials \( p_n(x;a,b,c,d|q) \) satisfying
\( \AW p_n(x;a,b,c,d|q) = c_n
p_{n-1}(x;aq^{1/2},bq^{1/2},cq^{1/2},dq^{1/2}|q) \), where
\( c_n = \tfrac{q^{n/2}-q^{-n/2}}{q^{1/2}-q^{-1/2}} \).

The forward shift operators for the orthogonal polynomials in the
Askey and \( q \)-Askey schemes are obtained from the divided difference
operator \( \D_{u,v} \) by the substitutions listed in
\Cref{tab:Duv-substitutions}. The corresponding orthogonal polynomial
families are listed in \Cref{tab:operator-families}.

Considering \( x=x(z)=(z+z^{-1})/2 \) as a function of \( z \) and
writing \( \check f(z) = f(x(z)) \), the Askey--Wilson operator may
also be written in the following form:
\[
  (\AW f)(x) =\frac{\check f(q^{1/2}z)-\check f(q^{-1/2}z)}
  {(q^{1/2}-q^{-1/2})(z-z^{-1})/2}.
\]
The Wilson, Racah, and \(q\)-Racah operators have analogous
representations. Considering \(x=x(s)\) as a function of \(s\) and
writing \(\check f(s)=f(x(s))\), we have
\begin{align}
\label{eq:wilson}  (\mathcal W f)(x)
  &=\frac{\check f(s+\ii/2)-\check f(s-\ii/2)}{2\ii s},
    && x=x(s)=s^2,\\
  (\mathcal R f)(x)
  &=\frac{\check f(s+1/2)-\check f(s-1/2)}
    {2s+\gamma+\delta+1},
    && x=x(s)=s(s+\gamma+\delta+1),\\
  (\Rq f)(x)
  &=\frac{\check f(s+1/2)-\check f(s-1/2)}
    {x(s+1/2)-x(s-1/2)},
    && x=x(s)=q^{-s}+\gamma\delta q^{s+1}.
\end{align}

\begin{table}
  \centering
  \renewcommand{\arraystretch}{1.4}
  \resizebox{\textwidth}{!}{%
  \begin{tabular}{c|cc|cc}
    operator & \( u(x) \) & \( v(x) \) & \(u(x)+v(x)\) & \(u(x)v(x)\) \\ \hline
    \( D \)
      & \( x \)
      & \( x \)
      & \(2x\)
      & \(x^2\) \\
    \( \Delta \)
      & \( x+1 \)
      & \( x \)
      & \(2x+1\)
      & \(x^2+x\) \\
    \( D_q \)
      & \( qx \)
      & \( x \)
      & \((q+1)x\)
      & \(qx^2\) \\
    \( \delta \)
      & \( x+\tfrac{\ii}{2} \)
      & \( x-\tfrac{\ii}{2} \)
      & \(2x\)
      & \(x^2+\tfrac14\) \\
    \( \mathcal{W} \)
      & \( x-\tfrac14+\ii\sqrt{x} \)
      & \( x-\tfrac14-\ii\sqrt{x} \)
      & \(2x-\tfrac12\)
      & \(x^2+\tfrac12x+\tfrac1{16}\) \\
    \( \mathcal{R} \)
      & \( x+\frac14+\sqrt{x+(\gamma+\delta+1)^2/4} \)
      & \( x+\frac14-\sqrt{x+(\gamma+\delta+1)^2/4} \)
      & \(2x+\tfrac12\)
      & \(x^2-\tfrac12x+\tfrac1{16}-\frac{(\gamma+\delta+1)^2}{4}\) \\
    \( \AW \)
      & \( \sigma x+\tau\sqrt{x^2-1} \)
      & \( \sigma x-\tau\sqrt{x^2-1} \)
      & \(2\sigma x\)
      & \(x^2+\tau^2\) \\
    \( \Rq \)
      & \( \sigma x+\tau\sqrt{x^2-4\gamma\delta q} \)
      & \( \sigma x-\tau\sqrt{x^2-4\gamma\delta q} \)
      & \(2\sigma x\)
      & \(x^2+4\gamma\delta q\tau^2\)
  \end{tabular}%
  }
  \caption{Substitutions realizing forward shift operators as
    \( \D_{u,v} \). Here, the \(\delta\) in the first column denotes
    the central difference operator, while the \(\delta\)'s in the
    other columns are parameters. The constants \( \sigma \) and
    \( \tau \) are defined by \( \sigma=\frac{q^{1/2}+q^{-1/2}}2 \)
    and \( \tau=\frac{q^{1/2}-q^{-1/2}}2 \).}
  \label{tab:Duv-substitutions}
\end{table}

\begin{table}
  \centering
  \renewcommand{\arraystretch}{1.15}
  \begin{tabular}{c|p{0.72\textwidth}}
    operator & families \\ \hline
    \(D\) & Hermite, Laguerre, Jacobi, and Bessel polynomials\\
    \(\Delta\) & Charlier, Meixner, Krawtchouk, and Hahn polynomials \\
    \(D_q\) & Al-Salam--Carlitz, \(q\)-Charlier, \(q\)-Meixner,
    \(q\)-Krawtchouk, \(q\)-Hahn, \(q\)-Laguerre, big and little \(q\)-Jacobi,
    and discrete \(q\)-Hermite polynomials \\
    \(\delta\) & continuous Hahn and Meixner--Pollaczek polynomials \\
    \(\mathcal{W}\) & Wilson and continuous dual Hahn polynomials \\
    \(\mathcal{R}\) & Racah and dual Hahn polynomials \\
    \(\AW\) & Askey--Wilson, Al-Salam--Chihara, continuous \(q\)-Hermite,
    continuous big \(q\)-Hermite,
    continuous \(q\)-Jacobi, and \(q\)-Meixner--Pollaczek polynomials \\
    \(\Rq\) & \(q\)-Racah, dual \(q\)-Hahn, and dual \(q\)-Krawtchouk polynomials
  \end{tabular}
  \caption{Orthogonal polynomial families associated with the forward shift operators.}
  \label{tab:operator-families}
\end{table}

Let \( T_n(x)= \check T_n(z) =(z^n+z^{-n})/2 \) and
\( U_{n-1}(x)=(z^n-z^{-n})/(z-z^{-1}) \) be the Chebyshev polynomials
of the first and second kind, where \( x=(z+z^{-1})/2 \). The
Askey--Wilson operator \( \AW \) and the corresponding averaging
operator \( \A_{\AW} \) behave well with the Chebyshev polynomials.

\begin{prop}\label{prop:AW-averaging-Chebyshev}
  We have
  \[
    \AW T_n(x)
    =
    \frac{q^{n/2}-q^{-n/2}}{q^{1/2}-q^{-1/2}}\,U_{n-1}(x),
    \qquad
    \A_{\AW}T_n(x)
    =\frac{q^{n/2}+q^{-n/2}}{2}\,T_n(x).
  \]
\end{prop}

\begin{proof}
  The first identity is \cite[(12.1.17)]{Ismail}.
  The second identity follows from the definition of
  \( \A_{\AW} \):
  \[
   \A_{\AW} T_n(x) = \frac{\check T_n(q^{1/2}z)+\check T_n(q^{-1/2}z)}{2} 
    = \frac{(q^{n/2}+q^{-n/2})T_n(x)}{2}. \qedhere
  \]
\end{proof}

\begin{remark}
  The formula for \( \A_{\AW}T_n \) is just the cosine addition
  formula. Indeed, write \( x=\cos\theta \) and
  \( q^{1/2}=e^{\ii\phi} \). Then
  \( \sigma=(q^{1/2}+q^{-1/2})/2=\cos\phi \) and
  \( \tau=(q^{1/2}-q^{-1/2})/2=\ii\sin\phi \). With the branch
  \( \sqrt{x^2-1}=\ii\sin\theta \), we have
  \[
    u(x)=\sigma x+\tau\sqrt{x^2-1}=\cos(\theta+\phi),
    \qquad
    v(x)=\sigma x-\tau\sqrt{x^2-1}=\cos(\theta-\phi).
  \]
  Thus \( \A_{\AW} \) averages the two shifts
  \( \theta\mapsto\theta+\phi \) and
  \( \theta\mapsto\theta-\phi \). Since
  \( T_n(\cos\theta)=\cos(n\theta) \), this gives
  \( \A_{\AW}T_n=((q^{n/2}+q^{-n/2})/2)T_n \).
\end{remark}

\subsection{Discriminants with two operators}

Recall from \Cref{def:disc} that given an operator \( T \), the
discriminant \( \disc(f;T) \) of a polynomial \( f \) of degree
\( n \) with leading coefficient \( A \) is defined by
\[
\disc(f;T)=(-1)^{\binom{n}{2}} A^{-1} \res(f,Tf).
\]
More generally, we define the discriminant associated with two
operators.

\begin{defn}\label{def:disc2}
  Let \( T_1 \) and \( T_2 \) be operators on polynomials. Suppose
  that \( f \) is a polynomial for which \( T_1 f \) has degree \( n \)
  with leading coefficient \( A \). We define
  the \emph{discriminant} \( \disc(f;T_1,T_2) \) of \( f \) with
  respect to \( (T_1,T_2) \) by
\[
\disc(f;T_1,T_2)=(-1)^{\binom{n}{2}} A^{-1} \res(T_1 f,T_2 f).
\]
\end{defn}

In this paper, we consider the two-operator discriminant in the case
\( T_1=\A_{u,v} \) and \( T_2=\D_{u,v} \). As in the case of the
ordinary discriminant, \( \disc(f;\A_{u,v},\D_{u,v}) \) is a
symmetric polynomial in the zeros of \( f \).

\begin{lem}\label{lem:disc2-degree}
  Let \( (u,v) \) be a degree-lowering pair with
  \[
    U(x) = u(x) + v(x) = a_1x+a_0, \qquad
    V(x) = u(x)v(x) = b_2x^2+b_1x+b_0.
  \]
  Let \(\theta_n\) be the leading coefficient of \(\A_{u,v}x^n\) and let
  \( f(x)=\prod_{i=1}^n(x-\alpha_i) \). Then
  \( \theta_n \disc(f;\A_{u,v},\D_{u,v})\) is a polynomial in the
  coefficients \( a_0,a_1,b_0,b_1,b_2 \) of \(U\) and \(V\) and the
  zeros \( \alpha_1,\ldots,\alpha_n \) of \( f \). In particular, for
  fixed \(u\) and \(v\), the discriminant
  \( \disc(f;\A_{u,v},\D_{u,v}) \) is a symmetric polynomial in the
  zeros of \(f\).
\end{lem}

\begin{proof}
  Let \( f(x)=\sum_{m=0}^n d_mx^m \). Since
  \(d_m=(-1)^{n-m}e_{n-m}(\alpha_1,\ldots,\alpha_n)\), where \(e_j\)
  denotes the \(j\)th elementary symmetric polynomial, every \(d_m\)
  is a symmetric polynomial in the zeros of \(f\).

  Let
  \[
    \A_{u,v}f=\sum_{k=0}^n p_kx^k,\qquad
    \D_{u,v}f=\sum_{k=0}^{n-1}q_kx^k .
  \]
  In the notation of \eqref{eq:HS}, we have
  \[
    \A_{u,v}f=\frac12\sum_{m=0}^n d_mS_m,\qquad
    \D_{u,v}f=\sum_{m=0}^n d_mH_m .
  \]
  The recurrence \eqref{eq:SH-rec} shows inductively that the
  coefficients of \(S_m\) and \(H_m\) are polynomials in
  \(a_0,a_1,b_0,b_1,b_2\). Hence \(p_k\) and \(q_k\) are polynomials in
  \(a_0,a_1,b_0,b_1,b_2,\alpha_1,\ldots,\alpha_n\), and are symmetric in
  \(\alpha_1,\ldots,\alpha_n\).
  The Sylvester determinant formula (see \cite[Section~14.6,~Exercise~29.(b)]{DummitFoote}) shows that
  \(\res(\A_{u,v}f,\D_{u,v}f)\) is a polynomial in the coefficients
  \(p_k\) and \(q_k\), hence has the same polynomiality and symmetry
  properties. By \Cref{def:disc2},
  \[
    \theta_n\disc(f;\A_{u,v},\D_{u,v})
    =
    (-1)^{\binom n2}\res(\A_{u,v}f,\D_{u,v}f)
  \]
  is therefore a polynomial in the five coefficients of \(U\) and \(V\)
  and the zeros of \(f\).
\end{proof}

Recall that the discriminant \( \disc(f;\D_u) \) has the product
formula in \Cref{pro:disc-prod} in terms of the zeros of \( f \).
There is an analogous product formula for
\( \disc(f;\A_{u,v},\D_{u,v}) \).

The authors received considerable help from ChatGPT 5.5 Pro in finding
and proving the following theorem. Note that this theorem implies that
if \( f \) is a polynomial of degree \( n \), then
\(\disc(f;\A_{u,v},\D_{u,v})\) is a symmetric polynomial in the
zeros of \( f \) of total degree \( n(n-1) \).

\begin{thm}\label{pro:disc-prod2}
  Let \( (u,v) \) be a degree-lowering pair with
  \( \A_{u,v}x = \rho x+\kappa \). Let \( \theta_n \) be the leading
  coefficient of \( \A_{u,v}x^n \). Then for a monic polynomial
  \( f(x)=\prod_{k=1}^n(x-\alpha_k) \), we have
  \begin{align}
    \disc(f;\A_{u,v},\D_{u,v})
    &=
    \frac{(2\rho)^{n(n-1)}}{\theta_n}
    \prod_{1\le i<j\le n}
    \left(
      u\left(\frac{\alpha_i+\alpha_j-2\kappa}{2\rho}\right)
      v\left(\frac{\alpha_i+\alpha_j-2\kappa}{2\rho}\right)
      -\alpha_i\alpha_j
    \right).
    \label{eq:avg-duv-disc-prod}
  \end{align}
\end{thm}

\begin{proof}
  By \Cref{lem:disc2-degree}, the left-hand side of
  \eqref{eq:avg-duv-disc-prod} multiplied by \( \theta_n \) is a
  polynomial in the coefficients \( a_0,a_1,b_0,b_1,b_2 \) of
  \( U=u+v \) and \( V=uv \) and in the zeros \( \alpha_i \); the same
  holds for the right-hand side multiplied by \( \theta_n \), since
  the factor \( (2\rho)^{n(n-1)}=a_1^{n(n-1)} \) clears the
  denominators \( (2\rho)^2 \) of the \( \binom n2 \) factors. We will
  prove the identity as polynomials in the \( \alpha_i \) and the
  coefficients of \( U \) and \( V \). Since two such polynomials that
  agree whenever \( b_2\ne0 \) are identical, we may assume
  \( \deg V=2 \), so that each polynomial \( G_k \) defined below is
  quadratic. Let
\[
  P(x)=(\A_{u,v}f)(x)=\frac{f(u(x))+f(v(x))}{2},\qquad
  Q(x)=(\D_{u,v}f)(x)=\frac{f(u(x))-f(v(x))}{u(x)-v(x)}.
\]
Then
\[
  f(u(x))=P(x)+\frac{u(x)-v(x)}2Q(x),\qquad
  f(v(x))=P(x)-\frac{u(x)-v(x)}2Q(x),
\]
and we have
\[
  P(x)^2-\frac{(u(x)-v(x))^2}{4}Q(x)^2
  =
  f(u(x))f(v(x))=   \prod_{k=1}^n G_k(x),
\]
where \( G_k(x)=(u(x)-\alpha_k)(v(x)-\alpha_k) \). Taking the
resultant with \(Q\), and using that adding a multiple of \(Q\) to the
first entry does not change the resultant, gives
\begin{equation}\label{eq:3}
  \res(P,Q)^2 = \prod_{k=1}^n\res(G_k,Q).
\end{equation}

We now compute \(\res(G_k,Q)\) for a fixed \(k\).  Let
\(\xi_1,\xi_2\) be the two roots of \(G_k\), and let \(c_k\) be the
leading coefficient of \(G_k\). Then
\begin{equation}\label{eq:8}
  G_k(x)=(u(x)-\alpha_k)(v(x)-\alpha_k) = c_k(x-\xi_1)(x-\xi_2).
\end{equation}
For \( \ell\in\{1,2\} \), since \( G_k(\xi_\ell)=0 \), we have
\( u(\xi_\ell)=\alpha_k \) or \( v(\xi_\ell)=\alpha_k \). We define
\( \beta_\ell \) by
\( \alpha_k+\beta_\ell = u(\xi_\ell)+ v(\xi_\ell) \). Then
\[
  Q(\xi_\ell)
  =\frac{f(u(\xi_\ell)) - f(v(\xi_\ell))}{u(\xi_\ell) - v(\xi_\ell)}
  =\frac{f(\beta_\ell) - f(\alpha_k)}{\beta_\ell-\alpha_k}
  =\frac{f(\beta_\ell)}{\beta_\ell-\alpha_k}.
\]
Since \(Q\) has degree \(n-1\), the formula for the resultant in terms
of the roots of \(G_k\), together with
\(f(\beta_\ell)=\prod_{j=1}^n(\beta_\ell-\alpha_j)\), gives
\begin{equation}\label{eq:4}
  \res(G_k,Q)
  =
  c_k^{\,n-1}\prod_{\ell=1}^2
  \frac{f(\beta_\ell)}{\beta_\ell-\alpha_k}
  =
  c_k^{\,n-1}
  \prod_{j\ne k}\prod_{\ell=1}^2(\beta_\ell-\alpha_j).
\end{equation}

On the other hand, since
\( \alpha_k+\beta_\ell = u(\xi_\ell)+ v(\xi_\ell) \) and
\(u(x)+v(x)=2\A_{u,v}x=2\rho x+2\kappa\), we have
\( \beta_\ell = 2\rho \xi_\ell +2\kappa -\alpha_k \).
Thus, letting \( x_{kj}=(\alpha_k+\alpha_j-2\kappa)/2\rho \),
we have
\[
  \prod_{\ell=1}^2(\beta_\ell-\alpha_j)
  = \prod_{\ell=1}^2
  \left(2\rho \xi_\ell +2\kappa -\alpha_k-\alpha_j\right)
  = (2\rho)^2 \prod_{\ell=1}^2 \left(\xi_\ell- x_{kj} \right)
  = \frac{(2\rho)^2}{c_k} G_k(x_{kj}).
\]
Since \(u(x_{kj})+v(x_{kj})=2\rho x_{kj}+2\kappa=\alpha_k+\alpha_j\),
by \eqref{eq:8}, we have
\[
  G_k(x_{kj})
  = u(x_{kj})v(x_{kj})-\alpha_k(u(x_{kj})+v(x_{kj})) + \alpha_k^2
  = u(x_{kj})v(x_{kj})-\alpha_k\alpha_j.
\]
Hence, we can rewrite \eqref{eq:4} as
\[
  \res(G_k,Q)
  =
  (2\rho)^{2(n-1)}
  \prod_{j\ne k} (u(x_{kj})v(x_{kj})-\alpha_k\alpha_j).
\]
This together with \eqref{eq:3} gives
\[
  \res(P,Q)^2
  =
  \left(
  (2\rho)^{n(n-1)}
  \prod_{1\le i<j\le n}
  \bigl(u(x_{ij})v(x_{ij})-\alpha_i\alpha_j\bigr)
  \right)^2.
\]
By \Cref{lem:disc2-degree,def:disc2}, \(\res(P,Q)\) is a polynomial in
the coefficients of \(U=u+v\) and \(V=uv\) and in the roots
\(\alpha_i\). Thus, \(\res(P,Q)\) is equal to the expression inside
the square brackets up to sign.

We determine this sign using the specialization \(u=v=x\). Then
\(\rho=1\), \(\kappa=0\), \(P=f\), and \(Q=f'\), so
  \[
    \res(P,Q)=\res(f,f')
    =
    (-1)^{\binom n2}
    \prod_{1\le i<j\le n}(\alpha_i-\alpha_j)^2.
  \]
  On the other hand, we have \(x_{ij}=(\alpha_i+\alpha_j)/2\), and
  hence
  \[
    (2\rho)^{n(n-1)}
  \prod_{1\le i<j\le n}
  \bigl(u(x_{ij})v(x_{ij})-\alpha_i\alpha_j\bigr)
    =
    \prod_{1\le i<j\le n}(\alpha_i-\alpha_j)^2.
  \]
  Comparing the two formulas gives
\[
  \res(\A_{u,v}f,\D_{u,v}f) = \res(P,Q) =
  (-1)^{\binom n2}(2\rho)^{n(n-1)}
  \prod_{1\le i<j\le n}
  \bigl(u(x_{ij})v(x_{ij})-\alpha_i\alpha_j\bigr).
\]
Finally, by \Cref{def:disc2}, we obtain \eqref{eq:avg-duv-disc-prod}.
\end{proof}

The following proposition shows that the discriminant
\( \disc(f;\A_{u,v},\D_{u,v}) \) generalizes \( \disc(f;\D_{u}) \) up
to scalar multiplication.

\begin{prop}\label{prop:linear-root-product}
  Let \(u(x)=ax+b\) and \(v(x)=cx+d\) be a degree-lowering pair, with
  \(c\ne0\), \(a+c\ne0\), and \(a^n+c^n\ne0\). Let
  \(f(x)=\prod_{i=1}^n(x-\alpha_i)\) be a monic polynomial. Then
  \begin{equation}\label{eq:5}
    \disc(f;\A_{u,v},\D_{u,v})
    =
    \frac{2c^{n(n-1)}}{a^n+c^n}
    \prod_{1\le i<j\le n}
    \bigl(\alpha_i-u(v^{-1}(\alpha_j))\bigr)
    \bigl(u(v^{-1}(\alpha_i))-\alpha_j\bigr).
  \end{equation}
  In particular, if \(v(x)=x\), then
  \[
    \disc(f;\A_{u,v},\D_{u,v})
    =\frac{2}{a^n+1}\disc(f;\D_{u}).
  \]
\end{prop}

\begin{proof}
  By \Cref{pro:disc-prod2} with \( 2\rho=a+c \), \( 2\kappa=b+d \),
  and \( \theta_n = (a^n+c^n)/2 \), we obtain
  \begin{align*}
    &\disc(f;\A_{u,v},\D_{u,v})\\
    &=
    \frac{2(a+c)^{n(n-1)}}{a^n+c^n}
    \prod_{1\leq i<j\leq n}
      \left(
	      \left( \frac{a(\alpha_i+\alpha_j-b-d)}{a+c}+b \right)
	      \left( \frac{c(\alpha_i+\alpha_j-b-d)}{a+c}+d \right)
      -\alpha_i\alpha_j
      \right)\\
    &=
    \frac{2}{a^n+c^n}
    \prod_{1\leq i<j\leq n}
    \bigl(c\alpha_i-a\alpha_j+ad-bc\bigr)
    \bigl(a\alpha_i-c\alpha_j+bc-ad\bigr).
  \end{align*}
  It is straightforward to check that this is equivalent
  to \eqref{eq:5}.
\end{proof}

The following corollary shows that the discriminants
\( \disc(f;\A_T, T) \) for the forward shift operators \( T \) in the
Askey and \( q \)-Askey schemes factor nicely.

\begin{cor}\label{cor:forward-shift-disc-prod}
Let \( f(x)=\prod_{k=1}^n(x-\alpha_k) \) be monic and put
\(N=\binom n2\).  For the forward shift operators in
\Cref{tab:Duv-substitutions}, the following product formulas hold,
where each product is over \(1\le i<j\le n\):
\begin{align*}
  \disc(f;\A_D,D)
  &= \prod_{i<j}(\alpha_i-\alpha_j)^2,\\
  \disc(f;\A_\Delta,\Delta)
  &= \prod_{i<j}
  (\alpha_i-\alpha_j-1)(\alpha_i-\alpha_j+1),\\
  \disc(f;\A_{D_q},D_q)
  &= \frac{2\prod_{i<j}
  (\alpha_i-q\alpha_j)(q\alpha_i-\alpha_j)}{1+q^n},\\
  \disc(f;\A_\delta,\delta)
  &= \prod_{i<j}
  (\alpha_i-\alpha_j+\ii)(\alpha_i-\alpha_j-\ii),\\
  \disc(f;\A_{\mathcal W},\mathcal W)
  &= \prod_{i<j}
  (s_i-s_j-\ii)(s_i-s_j+\ii)
  (s_i+s_j-\ii)(s_i+s_j+\ii),
  && \alpha_k=s_k^2,\\
  \disc(f;\A_{\mathcal R},\mathcal R)
  &= \prod_{i<j}
  (-s_i+s_j-1)(-s_i+s_j+1)
  (s_i+s_j-1)(s_i+s_j+1),
  && \alpha_k=s_k^2-\frac{(\gamma+\delta+1)^2}{4},\\
  \disc(f;\A_{\AW},\AW)
  &=
  \frac{2\prod_{i<j}
  (z_i-qz_j)(z_i-q^{-1}z_j)(z_iz_j-q)(z_iz_j-q^{-1})}
  {4^N(q^{n/2}+q^{-n/2})(z_1\cdots z_n)^{2(n-1)}},
  && \alpha_k=\frac{z_k+z_k^{-1}}2,\\
  \disc(f;\A_{\Rq},\Rq)
  &=
  \frac{2\prod_{i<j}
  (z_i-qz_j)(qz_i-z_j)(z_iz_j-qC)(qz_iz_j-C)}
  {q^{2N}(q^{n/2}+q^{-n/2})(z_1\cdots z_n)^{2(n-1)}},
  && C=\gamma\delta q,\quad \alpha_k=z_k+\frac{C}{z_k}.
\end{align*}
\end{cor}

\begin{proof}
  We will only give a proof for the Askey--Wilson operator since the
  other cases can be proved similarly. Let \( u \) and \( v \) be the
  functions with \( D_{u,v} = \AW \) given in
  \Cref{tab:Duv-substitutions}. By \Cref{pro:disc-prod2},
\[
  \disc(f;\A_{u,v},\D_{u,v})
  =
  \frac{(q^{1/2}+q^{-1/2})^{2N}}{(q^{n/2}+q^{-n/2})/2}
  \prod_{1\le i<j\le n}
  \left(
    \left(\frac{\alpha_i+\alpha_j}{q^{1/2}+q^{-1/2}}\right)^2
    +\left(\frac{q^{1/2}-q^{-1/2}}2\right)^2
    -\alpha_i\alpha_j
  \right).
\]
If \(\alpha_k=(z_k+z_k^{-1})/2\), we have
\[
  \left(\frac{\alpha_i+\alpha_j}{q^{1/2}+q^{-1/2}}\right)^2
    +\left(\frac{q^{1/2}-q^{-1/2}}2\right)^2
    -\alpha_i\alpha_j
    =
  \frac{
  (z_i-qz_j)(z_i-q^{-1}z_j)(z_iz_j-q)(z_iz_j-q^{-1})
  }
  {4(q^{1/2}+q^{-1/2})^2z_i^2z_j^2}.
\]
Hence, we obtain the desired formula.
\end{proof}

One may ask whether \(\disc(f;\D_{u,v})\) has a product formula
analogous to the one for
\(\disc(f;\A_{u,v},\D_{u,v})\) in \Cref{pro:disc-prod2}.
The following example suggests that this is unlikely.

\begin{exam}\label{exa:avg-duv-small-primes}
  Let \( f(x)=\prod_{i=1}^n(x-\alpha_i) \), and write
  \( U(x)=u(x)+v(x) \) and \( V(x)=u(x)v(x) \). We compute
  \(\disc(f;\A_{u,v},\D_{u,v})\) using \eqref{eq:HS}. For example, take
  \[
    \{\alpha_i\}=\{-2,0,3,7,11\},\qquad U(x)=-x,\qquad V(x)=-1.
  \]
  For this choice, \Cref{pro:disc-prod2} gives
  \[
    \disc(f;\A_{u,v},\D_{u,v})
    =-2\prod_{i<j}(1+\alpha_i\alpha_j).
  \]
  Thus its prime factors are controlled by simple pairwise expressions
  in the zeros. Indeed,
  \[
    \disc(f;\A_{u,v},\D_{u,v})
    =159279120
    =2^4\cdot3^2\cdot5\cdot7\cdot11\cdot13^2\cdot17.
  \]
  On the other hand, the discriminant \( \disc(f;\D_{u,v}) \)
  without the averaging operator \( \A_{u,v} \) is
  \[
    \disc(f;\D_{u,v})
    =40813644070955448
    =2^3\cdot3^3\cdot7^2\cdot11\cdot17\cdot31\cdot13103\cdot50767.
  \]
  In this example, the largest absolute value of a zero is \(11\),
  while the largest prime factor of the two-operator discriminant is
  \(17\). By contrast, without the averaging operator, the prime
  factors \(13103\) and \(50767\) occur. This contrast suggests that
  there is no simple product formula for \(\disc(f;\D_{u,v})\) in
  terms of the zeros of \( f \).
\end{exam}

\section{Linked pairs, resultants, and discriminants}
\label{sec:mixed-three-term}

In this section, we show that if a polynomial has a mixed three-term
relation with orthogonal polynomials, then their resultant has a
product formula. We show that if two OPS are related by a
divided difference operator \( \D_{u,v} \), then they satisfy a
mixed three-term relation, where the averaging operator \( \A_{u,v} \)
is applied to one OPS. We also show that if \( u \) and \( v \) are linear
functions, then there is also a mixed three-term relation without the
averaging operator.

\subsection{Mixed three-term relations}

The following is our key lemma, which shows that a mixed three-term
relation with an OPS gives a factorization of the corresponding
resultant.

\begin{lem}\label{lem:key}
  Let \( \{p_m(x)\} = \OPS(\{b_m\},\{\lambda_m\}) \) be a monic OPS.
  Fix \( n\ge2 \), and let \( f(x) \) be a polynomial of degree \( n \)
  satisfying
  \begin{equation}\label{eq:main}
    f(x) = C(x-A)p_{n-1}(x) - B p_{n-2}(x)
  \end{equation}
  for some constants \( A,B,C \). Then
  \begin{equation}\label{eq:main-resultant}
    \res(f,p_{n-1})
    = (-1)^{\binom{n}{2}} B^{n-1}\prod_{k=1}^{n-2}\lambda_k^k.
  \end{equation}
\end{lem}

\begin{proof}
  Let \( \beta_1,\dots,\beta_{n-1} \) be the zeros of \( p_{n-1}(x) \).
  Since \( p_{n-1}(x) \) is monic, \Cref{def:resultant} and
  \eqref{eq:main} give
  \[
    \res(f,p_{n-1})
    =(-1)^{n(n-1)}\prod_{k=1}^{n-1}f(\beta_k)
    =\prod_{k=1}^{n-1}\left(-Bp_{n-2}(\beta_k)\right)
    =(-B)^{n-1}\res(p_{n-1},p_{n-2}),
  \]
  where we used \( p_{n-1}(\beta_k)=0 \) and \( (-1)^{n(n-1)}=1 \).
  Applying Schur's lemma (\Cref{thm:schur}),
  we obtain \eqref{eq:main-resultant}.
\end{proof}

Now we introduce the notion of linked pairs, which plays an important role in this paper.

\begin{defn}\label{def:5}
  Let \( \{p_m(x)\} = \OPS(\{b_m\},\{\lambda_m\}) \) be a monic OPS.
  Let \( \{f_n(x)\} \) be a sequence of polynomials with
  \( \deg f_n(x) = n \). We say that \( (\{f_n(x)\}, \{p_n(x)\}) \) is
  a \emph{linked pair} if for all \( n\ge1 \),
  \[
    f_n(x) = C_n(x-A_n)p_{n-1}(x) - B_n p_{n-2}(x)
  \]
  for some constants \( A_n, B_n,C_n \),
  where \( p_{-1}(x)=0 \).
\end{defn}

The resultant of a linked pair has a simple product formula.

\begin{thm}\label{thm:linked-pair-resultant}
  Suppose that \( (\{f_n(x)\},\{p_n(x)\}) \) is a linked pair with
  \( \{p_m(x)\}=\OPS(\{b_m\},\{\lambda_m\}) \) and
  \[
    f_n(x)=C_n(x-A_n)p_{n-1}(x)-B_np_{n-2}(x)
    \qquad (n\ge1),
  \]
  where \(p_{-1}(x)=0\). Then, for all \(n\ge2\),
  \[
    \res(f_n,p_{n-1})
    =
    (-1)^{\binom n2}B_n^{\,n-1}
    \prod_{k=1}^{n-2}\lambda_k^k.
  \]
\end{thm}

\begin{proof}
  This follows immediately from \Cref{lem:key} applied to the
  defining relation of the linked pair with \(C=C_n\), \(A=A_n\), and
  \(B=B_n\).
\end{proof}

We introduce another useful definition for pairs of OPS.

\begin{defn}\label{def:3}
  Let \( T \) be a degree-lowering operator and let \( \eta_n \) be
  the leading coefficient of \( T x^n \). For monic OPS
  \( \{f_n(x)\} \) and \( \{\widetilde{f}_n(x)\} \), we say that
  \( (\{f_n(x)\}, \{\widetilde{f}_n(x)\}) \) is a
  \emph{\( T \)-pair} if
 \[
   T f_n = \eta_n\widetilde{f}_{n-1} \qquad (n\ge1).
 \]
 In this case, we will also say that \( \{f_n\} \) and \( \{\wf_n\} \)
 form a \( T \)-pair.
 If \( T=\D_{u,v} \), we call this a \( \D_{u,v} \)-pair.
\end{defn}

Note that each OPS in the Askey and \( q \)-Askey schemes gives rise
to a \( T \)-pair, where \( T \) is its forward shift operator.
For example,
if 
\( f_n(x) = p_n(x;a,b,c,d|q) \) is the Askey--Wilson polynomial
and 
\[
  \widetilde{f}_n(x)=
  p_n(x;aq^{1/2},bq^{1/2},cq^{1/2},dq^{1/2}|q),
\]
then \( (\{f_n(x)\}, \{\widetilde{f}_n(x)\}) \) is a
\( \AW \)-pair.

The following is the main result in this paper, which shows that a
pair of OPS is a \( \D_{u,v} \)-pair if and only if they satisfy a
mixed three-term relation. Recall that if \( (u,v) \) is a
degree-lowering pair, then \( \mathcal{A}_{u,v}x = \rho x+\kappa \)
for some constants \( \rho \) and \( \kappa \).

\begin{thm}[Main Result]\label{thm:mixed-3-term}
  Let \( (u,v) \) be a degree-lowering pair such that
  \( \mathcal{A}_{u,v}x = \rho x+\kappa \), and let \( \eta_n \)
  denote the leading coefficient of \( \D_{u,v} x^n \). Let
  \[
    \{f_n(x)\} = \OPS(\{b_n\},\{\lambda_n\}),\qquad 
    \{\widetilde{f}_n(x)\} =
    \OPS(\{\widetilde{b}_n\},\{\widetilde{\lambda}_n\}).
  \]
  Then \( (\{f_n(x)\}, \{\widetilde{f}_n(x)\}) \) is a
  \( \D_{u,v} \)-pair if and only if, for all \( n \geq 1 \),
  \begin{equation}\label{eq:Af}
   \A_{u,v} f_n
   = ((\eta_{n+1}-\rho\eta_n)x
     -\eta_{n+1}\widetilde{b}_{n-1}
     +\eta_n b_n-\kappa\eta_n)\wf_{n-1} 
     -(\eta_{n+1}\widetilde{\lambda}_{n-1}
     -\eta_{n-1}\lambda_n)\wf_{n-2},
 \end{equation}
 where we set
 \( \eta_0 = \widetilde{f}_{-1} = \widetilde{\lambda}_0 = 0 \). In
 particular, if \( (\{f_n(x)\}, \{\widetilde{f}_n(x)\}) \) is a
 \( \D_{u,v} \)-pair, then
 \( (\{\A_{u,v}f_n(x)\}, \{\widetilde{f}_n(x)\}) \) is a linked pair.
\end{thm}

\begin{proof}
  Suppose first that \( (\{f_n(x)\},\{\widetilde f_n(x)\}) \) is a
  \( \D_{u,v} \)-pair. Applying \( \D_{u,v} \) to
  \( f_{n+1} = (x-b_n)f_n - \lambda_n f_{n-1} \) via the product
  rule~\eqref{eq:avg-product}, and using \( \D_{u,v}x = 1 \),
  \( \mathcal{A}_{u,v}(x-b_n) = \rho x + \kappa - b_n \), and
  \( \D_{u,v}f_n = \eta_n\widetilde{f}_{n-1} \), we obtain, for
  \( n \geq 1 \),
  \[
    \eta_{n+1}\widetilde{f}_n
    = \eta_n(\rho x + \kappa - b_n)\widetilde{f}_{n-1}
      + \mathcal{A}_{u,v}f_n - \eta_{n-1}\lambda_n\widetilde{f}_{n-2}.
  \]
  where for \( n=1 \) we use
  \(\D_{u,v}f_0=0 = \eta_0 \widetilde{f}_{-1}\). Solving for
  \( \mathcal{A}_{u,v}f_n \) and substituting the three-term
  recurrence for \( \widetilde{f}_n \) gives~\eqref{eq:Af}.

  Conversely, suppose that~\eqref{eq:Af} holds for all \( n \geq 1 \).
  We prove by induction on \(n\) that
  \( \D_{u,v}f_n=\eta_n\wf_{n-1} \). For \(n=1\), this is immediate
  since \(f_1=x-b_0\), \( \D_{u,v}x=1 \), \( \eta_1=1 \), and
  \( \wf_0=1 \). Now let \(n\ge1\) and assume that
  \[
    \D_{u,v}f_n=\eta_n\wf_{n-1}, \qquad
    \D_{u,v}f_{n-1}=\eta_{n-1}\wf_{n-2},
  \]
  where for \( n=1 \) the second identity is understood as
  \( \D_{u,v}f_0 = 0 = \eta_0 \widetilde{f}_{-1} \). Applying
  \( \D_{u,v} \) to \( f_{n+1}=(x-b_n)f_n-\lambda_n f_{n-1} \) and
  then using \eqref{eq:Af}, we obtain
  \[
    \D_{u,v}f_{n+1}
    =\eta_{n+1}(x-\widetilde b_{n-1})\wf_{n-1}
      -\eta_{n+1}\widetilde\lambda_{n-1}\wf_{n-2}
    =\eta_{n+1}\wf_n.
  \]
  Hence \( (\{f_n(x)\}, \{\widetilde{f}_n(x)\}) \) is a
  \( \D_{u,v} \)-pair.
\end{proof}

For \(n\ge2\), since the polynomials \(x\wf_{n-1}\), \( \wf_{n-1} \), and
\( \wf_{n-2} \) have degrees \(n\), \(n-1\), and \(n-2\),
respectively, the coefficients of the mixed three-term relation
\eqref{eq:Af} are uniquely determined if they exist.
For \(n=1\), the coefficient of \(\wf_{-1}=0\) is irrelevant.

As a corollary of \Cref{thm:mixed-3-term}, if
\( (\{f_n\}, \{\wf_n\}) \) is a \( \D_{u,v} \)-pair, we obtain an
explicit product formula for the discriminant
\( \disc(f_n;\mathcal{A}_{u,v},\D_{u,v}) \).

\begin{cor}\label{cor:avg-T-disc}
  Let \(n\ge2\), and let \( (u,v) \) be a degree-lowering pair such that
  \( \mathcal{A}_{u,v}x = \rho x+\kappa \), and let \( \eta_n \)
  denote the leading coefficient of \( \D_{u,v} x^n \). Suppose that
  \[
    \{f_n(x)\} = \OPS(\{b_n\},\{\lambda_n\}),\qquad 
    \{\widetilde{f}_n(x)\} =
    \OPS(\{\widetilde{b}_n\},\{\widetilde{\lambda}_n\})
  \]
  form a \( \D_{u,v} \)-pair and
  \( \eta_{n+1}-\rho \eta_n\ne0 \). Then
  \begin{equation}\label{eq:avg-T-disc}
    \disc(f_n;\mathcal{A}_{u,v},\D_{u,v})
    = \frac{\eta_n^{\,n} (\eta_{n+1}\widetilde{\lambda}_{n-1} - \eta_{n-1}\lambda_n)^{\,n-1}
      \prod_{k=1}^{n-2}\widetilde{\lambda}_k^{\,k}}
      {\eta_{n+1}-\rho \eta_n}.
  \end{equation}
\end{cor}

\begin{proof}
  This follows from \Cref{thm:mixed-3-term,lem:key} and the definition
  of \( \disc(f_n;T_1,T_2) \) in \Cref{def:disc2}, using that the leading
  coefficient of \( \mathcal{A}_{u,v}f_n \) in \eqref{eq:Af} is
  \( \eta_{n+1}-\rho \eta_n \).
\end{proof}

\begin{remark}\label{rem:1}
Our method is different from Ismail's method, which
finds the following relation (see \cite[Theorem~3.2.1]{Ismail})
\begin{equation}\label{eq:13}
Tp_N(x)=-B_N(x) p_N(x)+A_N(x) p_{N-1}(x),
\end{equation}
where \( B_N(x) \) and \( A_N(x) \) are rational functions in \( x \),
and apply Schur's lemma \cite[Lemma~3.4.1]{Ismail}. This method works
if \( \prod_{i=1}^N A_N(x_{N,i}) \) can be evaluated, where
\(x_{N,1},\ldots,x_{N,N}\) are the zeros of \(p_N\). Finding
\( A_N(x) \) and \( B_N(x) \) requires evaluating integrals or
summations; see \cite[Theorem 3.2.1]{Ismail} for \( T=d/dx \),
\cite[Theorem 6.3.1]{Ismail} for \( T=\Delta \), and \cite[Theorem
18.6.1]{Ismail} for \( T=\AW \).
\end{remark}

\subsection{The linear case}

An important special case of a degree-lowering pair \( (u,v) \) occurs
when both \( u \) and \( v \) are linear polynomials. For our
purposes, it suffices to consider \( u(x)=\alpha x+\beta \) and
\( v(x)=x+\gamma \). We will show that, in this case, the mixed
three-term relation for \( \mathcal{A}_{u,v} f_n \) in
\Cref{thm:mixed-3-term} can be rewritten as a relation for \( f_n \).

Let \( [n]_\alpha = 1+\alpha+ \cdots +\alpha^{n-1} \).

\begin{lem}\label{lem:Duv-monomial}
  Let \( u(x)=\alpha x+\beta \) and \( v(x)=x+\gamma \). Then, for
  \( n\ge1 \), the coefficient of \(x^{n-1}\) in
  \( \D_{u,v}(x^n) \) is \( [n]_\alpha \).
\end{lem}

\begin{proof}
  The definition of \( \D_{u,v} \) gives
  \begin{equation}\label{eq:D_uv_xn}
    \D_{u,v}(x^n)=\left.\frac{y^n-z^n}{y-z} \right|_{y=u,z=v}
    =\sum_{i=0}^{n-1}u(x)^i v(x)^{n-1-i}.
  \end{equation}
  The coefficient of \(x^{n-1}\) in \( u(x)^i v(x)^{n-1-i} \) is
  \(\alpha^i\). Summing over \(i\) gives
  \( \sum_{i=0}^{n-1}\alpha^i=[n]_\alpha \).
\end{proof}

In particular, if \((u,v)\) is a degree-lowering pair, then
\([n]_\alpha\ne0\) and it is the leading coefficient of
\(\D_{u,v}(x^n)\).

Note that if \( v(x) = x+\gamma \) and
\( \{\widetilde{f}_n(x)\} =
\OPS(\{\widetilde{b}_n\},\{\widetilde{\lambda}_n\}) \), then
considering the three-term recurrence of \( \widetilde{f}_n(x) \), we
have
\[
  \{(\widetilde f_n \circ v^{-1})(x)\} =
  \{\widetilde f_n(x-\gamma)\} =\OPS(\{\widetilde b_n+\gamma\},
  \{\widetilde\lambda_n\}).
\]
The following theorem gives a mixed three-term relation for
\( f_n \) when \( u(x)=\alpha x+\beta \) and \( v(x)=x+\gamma \).

\begin{thm}\label{lem:linear-companion}
  Let \( u(x)=\alpha x+\beta \) and \( v(x)=x+\gamma \), and let
  \( \eta_n \) denote the leading coefficient of \( \D_{u,v} x^n \).
  Suppose that
  \[
    \{f_n(x)\} = \OPS(\{b_n\},\{\lambda_n\}),\qquad 
    \{\widetilde{f}_n(x)\} =
    \OPS(\{\widetilde{b}_n\},\{\widetilde{\lambda}_n\})
  \]
  form a \( \D_{u,v} \)-pair. Let
  \[
    \{g_n(x)\} = \{\widetilde f_n(x-\gamma)\}
    = \OPS(\{\widetilde b_n+\gamma\}, \{\widetilde\lambda_n\}).
  \]
  Then, for all \( n \geq 1 \),
  \begin{equation}\label{eq:linear-conn}
    f_n
    = \left(x-(\eta_{n+1}\widetilde b_{n-1}-\eta_nb_n+\beta \eta_n+\gamma)\right)g_{n-1}
      - \left(\eta_{n+1}\widetilde\lambda_{n-1}-\eta_{n-1}\lambda_n\right)g_{n-2},
  \end{equation}
  where we set \( \eta_0 = g_{-1} = \widetilde{\lambda}_0 = 0 \).
\end{thm}

\begin{proof}
  By \Cref{lem:Duv-monomial}, \( \eta_n=[n]_\alpha \), and hence
  \( \eta_{n+1}=\alpha\eta_n+1 \). Applying the non-symmetric product
  rule~\eqref{eq:leibniz} to
  \( f_{n+1}=(x-b_n)f_n-\lambda_n f_{n-1} \), for \(n\ge1\), gives
  \[
    \eta_{n+1}\widetilde f_n
    = (u-b_n)\eta_n\widetilde f_{n-1}
      + f_n\circ v
      - \lambda_n\eta_{n-1}\widetilde f_{n-2},
  \]
  where for \( n=1 \) we use
  \(\D_{u,v}f_0=0=\eta_0 \widetilde{f}_{-1}\). Thus, using
  \( \widetilde f_n=(x-\widetilde b_{n-1})\widetilde f_{n-1}
  -\widetilde\lambda_{n-1}\widetilde f_{n-2} \), we get
  \begin{align*}
    f_n\circ v
    &= \eta_{n+1}\widetilde f_n
      - \eta_n(u-b_n)\widetilde f_{n-1}
      + \eta_{n-1}\lambda_n\widetilde f_{n-2}\\
    &= \left(\eta_{n+1}\,x - \eta_n\,u
      - \eta_{n+1}\widetilde b_{n-1}
      + \eta_n b_n\right)\widetilde f_{n-1}
      - \left(\eta_{n+1}\widetilde\lambda_{n-1}
      - \eta_{n-1}\lambda_n\right)\widetilde f_{n-2}.
  \end{align*}
  Replacing \( x \) by \( v^{-1}(x) \) gives
  \[
    f_n
    = \left(\eta_{n+1} v^{-1}(x) - \eta_n\,u \circ v^{-1}(x)
      - \eta_{n+1}\widetilde b_{n-1}
      + \eta_n b_n\right) g_{n-1}
      - \left(\eta_{n+1}\widetilde\lambda_{n-1}
      - \eta_{n-1}\lambda_n\right) g_{n-2}.
  \]
  Since \( v^{-1}(x)=x-\gamma \),
  \( u\circ v^{-1}(x)=\alpha(x-\gamma)+\beta \), and
  \( \eta_{n+1}-\alpha\eta_n=1 \), the coefficient of \(g_{n-1}\) is
  \[
    \eta_{n+1}(x-\gamma) - \eta_n\left(\alpha(x-\gamma)+\beta\right)
    - \eta_{n+1}\widetilde b_{n-1} + \eta_n b_n
    = x-(\eta_{n+1}\widetilde b_{n-1}-\eta_nb_n
    +\beta \eta_n+\gamma),
  \]
  which gives~\eqref{eq:linear-conn}.
\end{proof}

The following corollary shows that under certain conditions,
\( \disc(f_n;\mathcal{A}_{u,v},\D_{u,v}) \) can be written as a
one-operator discriminant.

\begin{cor}\label{cor:linear-disc-ratio}
  Let \( u(x)=\alpha x+\beta \) and \( v(x)=x+\gamma \), and let
  \( \eta_n \) denote the leading coefficient of \( \D_{u,v} x^n \).
  Suppose that
  \[
    \{f_n(x)\} = \OPS(\{b_n\},\{\lambda_n\}),\qquad 
    \{\widetilde{f}_n(x)\} =
    \OPS(\{\widetilde{b}_n\},\{\widetilde{\lambda}_n\})
  \]
  form a \( \D_{u,v} \)-pair. Let
  \[
    \{g_n(x)\}=\{\widetilde f_n(x-\gamma)\}
    =\OPS(\{\widetilde b_n+\gamma\}, \{\widetilde\lambda_n\}),
  \]
  and let \( S \) be the operator defined by
  \( S f= (\D_{u,v} f)\circ v^{-1} \). Then, for \(n\ge2\),
  \[
    \disc(f_n;S)
    =
    \eta_n^{\,n}
    \left(\eta_{n+1}\widetilde\lambda_{n-1}
      -\eta_{n-1}\lambda_n\right)^{n-1}
    \prod_{k=1}^{n-2}\widetilde\lambda_k^{\,k}.
  \]
  In particular, if \( \eta_{n+1}-\rho \eta_n\ne0 \),
  where \( \rho=(\alpha+1)/2 \), then
  \[
    \disc(f_n;\mathcal{A}_{u,v},\D_{u,v})
    = \frac{1}{\eta_{n+1}-\rho \eta_n}\disc(f_n;S).
  \]
\end{cor}

\begin{proof}
  By \Cref{lem:linear-companion}, the sequence
  \( \{g_n(x)\}=\OPS(\{\widetilde b_n+\gamma\}, \{\widetilde\lambda_n\}) \)
  is a monic OPS, and \( Sf_n=\eta_ng_{n-1} \). Applying \Cref{lem:key}
  to~\eqref{eq:linear-conn} gives
  \[
    \disc(f_n;S)
    = \eta_n^{\,n}
      (\eta_{n+1}\widetilde\lambda_{n-1}-\eta_{n-1}\lambda_n)^{n-1}
      \prod_{k=1}^{n-2}\widetilde\lambda_k^{\,k}.
  \]
  The ratio follows by comparing this formula with \Cref{cor:avg-T-disc}.
\end{proof}

\subsection{Applications to Askey schemes}
\label{sec:appl-askey-schem}

In this subsection, we apply the results in the previous subsection to
the orthogonal polynomials \( p_n(x) \) in the Askey and \( q \)-Askey
schemes. Since each forward shift operator in the Askey and
\( q \)-Askey schemes can be realized as \( \D_{u,v} \), we obtain a
product formula for the two-operator discriminant
\( \disc(p_n(x);\A_{u,v},\D_{u,v}) \). If
\(u(x)=\alpha x+\beta\) and \(v(x)=x+\gamma\), we also obtain a product
formula for the one-operator discriminant \(\disc(p_n(x);S)\), where
\(Sf=(\D_{u,v}f)\circ v^{-1}\).

\begin{table}
  \centering
  \renewcommand{\arraystretch}{1.6}
  \setlength{\tabcolsep}{2pt}
  \begin{tabular}{l|cccccc}
    \( \D_{u,v} \) & \( u(x) \) & \( v(x) \) & \( \rho \) & \( \kappa \) & \( \eta_n \)
    & \( \eta_{n+1}-\rho \eta_n \) \\ \hline
    \( D \) & \( x \) & \( x \) & \( 1 \) & \( 0 \) & \( n \) & \( 1 \) \\
    \( \Delta \) & \( x+1 \) & \( x \) & \( 1 \) & \( \tfrac12 \) & \( n \) & \( 1 \) \\
    \( D_q \) & \( qx \) & \( x \) & \( \tfrac{q+1}{2} \) & \( 0 \) & \( [n]_q \)
      & \( \tfrac{1+q^n}{2} \) \\
    \( \delta \) & \( x+\frac{\ii}{2} \) & \( x-\frac{\ii}{2} \) & \( 1 \) & \( 0 \) & \( n \) & \( 1 \) \\
    \( \mathcal{W} \) & nonlinear & nonlinear & \( 1 \) & \( -\tfrac14 \) & \( n \) & \( 1 \) \\
    \( \mathcal{R} \) & nonlinear & nonlinear & \( 1 \) & \( \tfrac14 \) & \( n \) & \( 1 \) \\
    \( \AW \) & nonlinear & nonlinear & \( \tfrac{q^{1/2}+q^{-1/2}}{2} \) & \( 0 \) &
      \( \tfrac{q^{n/2}-q^{-n/2}}{q^{1/2}-q^{-1/2}} \) &
      \( \tfrac{q^{n/2}+q^{-n/2}}{2} \) \\
    \( \Rq \) & nonlinear & nonlinear & \( \tfrac{q^{1/2}+q^{-1/2}}{2} \) & \( 0 \) &
      \( \tfrac{q^{n/2}-q^{-n/2}}{q^{1/2}-q^{-1/2}} \) &
      \( \tfrac{q^{n/2}+q^{-n/2}}{2} \)
  \end{tabular}
  \caption{The functions \( u(x) \) and \( v(x) \), and the constants \( \rho \),
  \( \kappa \), \( \eta_n \), and \( \eta_{n+1}-\rho \eta_n \) for the
  realizations of the Askey-scheme operators as \( \D_{u,v} \).
  Here, \( \A_{u,v}x = \rho x+\kappa \)
and \( \eta_n \) is the leading coefficient of \( \D_{u,v} x^n \).}
  \label{tab:askey-duv-constants}
\end{table}

\begin{lem}\label{lem:table-eta}
  Let \( \D_{u,v} \) be a forward shift operator in the Askey or
  \( q \)-Askey scheme with \( \A_{u,v}x = \rho x +\kappa \). Let
  \( \eta_n \) be the leading coefficients of \( \D_{u,v}x^n \). Then
  the values of \( \eta_n \), \( \kappa \), and
  \( \eta_{n+1}-\rho \eta_n \) are those listed in
  \Cref{tab:askey-duv-constants}.
\end{lem}

\begin{proof}
  Since the polynomials \( f_n \) and \( \widetilde f_{n-1} \) are
  monic in the relation \( \D_{u,v}f_n=\eta_n\widetilde f_{n-1} \),
  the constant \( \eta_n \) is the leading coefficient of
  \( \D_{u,v}x^n \). Recall~\eqref{eq:D_uv_xn} that
  \[
    \D_{u,v}x^n
    =\sum_{j=0}^{n-1}u(x)^jv(x)^{n-1-j}.
  \]
  Thus, if the leading linear terms of \( u \) and \( v \) are
  \( a x \) and \( b x \), respectively, then
  \( \eta_n=\sum_{j=0}^{n-1}a^jb^{n-1-j} \). For \( D \), the limiting
  definition gives \( \eta_n=n \). For \( \Delta \), \( \delta \),
  \( \mathcal W \), and \( \mathcal R \), we have \( a=b=1 \), so
  \( \eta_n=n \). For \( D_q \), we have \( a=q \) and \( b=1 \), so
  \( \eta_n=[n]_q \). Finally, for \( \AW \) and
  \( \Rq \), the leading linear coefficients are \( a=q^{1/2} \) and
  \( b=q^{-1/2} \), and hence
  \[
    \eta_n=\sum_{j=0}^{n-1}q^{j/2}q^{-(n-1-j)/2}
    =\frac{q^{n/2}-q^{-n/2}}{q^{1/2}-q^{-1/2}}.
  \]
  The entries in the final column follow by substituting these formulas
  for \( \eta_n \).  
\end{proof}

For any OPS \( \{p_n(x)\} \) in the Askey and \( q \)-Askey schemes,
we can evaluate the discriminant \( \disc(p_n(x);\A_T, T) \) using
\Cref{cor:avg-T-disc}. If \(T=\D_{u,v}\) with \(u(x)=\alpha x+\beta\)
and \(v(x)=x+\gamma\), we can also evaluate \(\disc(p_n(x);S)\), where
\(Sf=(Tf)\circ v^{-1}\), using \Cref{cor:linear-disc-ratio}. To do
this, we need to compute \( \eta_n \), \( \eta_{n+1}-\rho \eta_n \),
and \( \eta_{n+1}\widetilde\lambda_{n-1} -\eta_{n-1}\lambda_n \). In
\Cref{sec:app}, we list the forward shift operator \( T \), and the
values of \( \widetilde\lambda_n \) and
\( \eta_{n+1}\widetilde\lambda_{n-1} -\eta_{n-1}\lambda_n \) for each
OPS. The values of \( \eta_n \) and \( \eta_{n+1}-\rho \eta_n \) are
given in \Cref{tab:askey-duv-constants}.

\begin{exam}[Askey--Wilson]\label{ex:askey-wilson-disc}
  Let \( f_n = p_n(x;a,b,c,d\mid q)\) be the monic Askey--Wilson
  polynomial. The Askey--Wilson operator satisfies
  \( \AW f_n = \eta_n \widetilde{f}_{n-1} \) with
  \( \eta_n=(q^{n/2}-q^{-n/2})/(q^{1/2}-q^{-1/2}) \) and
  \[
    \{\widetilde f_n\}
    =\{p_n(x;aq^{1/2},bq^{1/2},cq^{1/2},dq^{1/2}\mid q)\}
    =\OPS(\{\widetilde b_n\}, \{\widetilde\lambda_n\}),
  \]
  where \( \widetilde b_n \) and \( \widetilde\lambda_n \) are obtained
  from \eqref{eq:AW-b,la} by replacing \( (a,b,c,d) \) with
  \( (aq^{1/2},bq^{1/2},cq^{1/2},dq^{1/2}) \). For \(n\geq2\), by
  \Cref{tab:askey-duv-constants} and \Cref{cor:avg-T-disc}, we obtain
  \[
    \disc(p_n(x;a,b,c,d\mid q);\A_{\AW},\AW)
    =
    \frac{2\eta_n^n}{q^{n/2}+q^{-n/2}}\,
    B_n^{n-1}\prod_{m=1}^{n-2}\widetilde{\lambda}_m^m,
  \]
  where
  \begin{align*}
    B_n
    &=\eta_{n+1}\widetilde{\lambda}_{n-1}
      -\eta_{n-1}\lambda_n\\
    &=
    \frac{q^{-n/2}(1-q^{n-1})(1-abq^{n-1})
    (1-acq^{n-1})(1-adq^{n-1})(1-bcq^{n-1})(1-bdq^{n-1})(1-cdq^{n-1})}
    {4(1-abcdq^{2n-3})(1-abcdq^{2n-2})^2} ,
  \end{align*}
  and
  \begin{align*}
    \widetilde{\lambda}_m
    &=
    \frac{(1-q^m)(1-abq^m)(1-acq^m)(1-adq^m) (1-abcdq^m)(1-bcq^m)(1-bdq^m)(1-cdq^m)}
    {4(1-abcdq^{2m-1})(1-abcdq^{2m})^2(1-abcdq^{2m+1})} .
  \end{align*}
\end{exam}

\begin{exam}[Continuous big \(q\)-Hermite]
  \label{ex:continuous-big-q-hermite-disc}
  Let \( f_n=H_n(x;a\mid q) \) be the monic continuous big
  \(q\)-Hermite polynomial. The Askey--Wilson operator satisfies
  \begin{equation}\label{eq:big_q_hermite}
  \AW f_n = \eta_n \wf_{n-1}.
  \end{equation}
  where
  \( \eta_n=(q^{n/2}-q^{-n/2})/(q^{1/2}-q^{-1/2}) \) and
  \( \wf_{n} = H_{n}(x;aq^{1/2}\mid q) \).
  We have \( \{f_n\} =\OPS(\{b_n\}, \{\lambda_n\}) \)
  and \( \{\widetilde{f}_n\} =\OPS(\{\widetilde{b}_n\}, \{\widetilde{\lambda}_n\}) \), where
  \[
    b_n = \frac{aq^n}{2}, \qquad 
    \widetilde{b}_n = \frac{aq^{n+1/2}}{2}, \qquad 
    \lambda_n=\widetilde\lambda_n=\frac{1-q^n}{4}.
  \]
  Hence, for \(n\geq2\), by
  \Cref{tab:askey-duv-constants} and \Cref{cor:avg-T-disc}, we obtain
  \[
    \disc(H_n(x;a\mid q);\A_{\AW},\AW)
    =
    \frac{2\eta_n^n}{q^{n/2}+q^{-n/2}}
    \left(\frac{q^{-n/2}(1-q^{n-1})}{4}\right)^{n-1}
    \prod_{m=1}^{n-2}\left(\frac{1-q^m}{4}\right)^m.
  \]

  We note that the following resultant with the squared Askey--Wilson
  operator also factors:
  \[
    \res\!\left(H_n(x;a\mid q),\AW^2H_{n+1}(x;a\mid q)\right)
    =
    (-1)^{\binom n2}(\eta_{n+1}\eta_n)^n
    \prod_{m=1}^{n-1}\left(\frac{1-q^m}{4}\right)^m.
  \]
  By \eqref{eq:big_q_hermite}, this is equivalent to
  \[
    \res\!\left(H_n(x;a\mid q),H_{n-1}(x;aq\mid q)\right)
    =
    (-1)^{\binom n2}
    \prod_{m=1}^{n-1}\left(\frac{1-q^m}{4}\right)^m,
  \]
  which follows from \Cref{linearthm2} in \Cref{sec:result-orth-polyn},
  applied with \( A=(a+a^{-1})/2 \), using
  \( H_n(A;a\mid q)=(2a)^{-n} \).
\end{exam}

Finally, we give an example using \Cref{cor:linear-disc-ratio}.

\begin{exam}[Meixner--Pollaczek]\label{ex:meixner-pollaczek-disc}
  Let \( f_n = P_n^{(\lambda)}(x;\phi)\) be the monic
  Meixner--Pollaczek polynomial. We have
  \( \delta f_n =n \widetilde{f}_{n-1} \), where \( \delta \) is the
  operator \( \D_{u,v} \) with \(u(x)=x+\ii/2\) and \(v(x)=x-\ii/2\),
  and \( \widetilde f_n=P_n^{(\lambda+1/2)}(x;\phi) \).
  We have \( \{f_n\} =\OPS(\{b_n\}, \{\lambda_n\}) \) and
  \( \{\widetilde{f}_n\}
  =\OPS(\{\widetilde{b}_n\}, \{\widetilde{\lambda}_n\}) \), where
  \begin{align*}
    b_n&=-\frac{(n+\lambda)\cos\phi}{\sin\phi},
    &\lambda_n&=\frac{n(n+2\lambda-1)}{4\sin^2\phi},\\
    \widetilde b_n&=-\frac{(n+\lambda+1/2)\cos\phi}{\sin\phi},
    &\widetilde\lambda_n&=\frac{n(n+2\lambda)}{4\sin^2\phi}.
  \end{align*}
  Let \( S f=(\delta f)\circ v^{-1} \). Then, for \(n\geq2\),
  \Cref{cor:linear-disc-ratio} gives
  \[
    \disc(P_n^{(\lambda)};S)
    =
    n^n B_n^{n-1}
    \prod_{m=1}^{n-2}
    \left(\frac{m(m+2\lambda)}{4\sin^2\phi}\right)^m,
  \]
  where
  \[
    B_n
    =\eta_{n+1}\widetilde\lambda_{n-1}
    -\eta_{n-1}\lambda_n
    =
    \frac{(n-1)(n+2\lambda-1)}{4\sin^2\phi}.
  \]
  Since \( \eta_{n+1}-\rho\eta_n=1 \), by
  \Cref{cor:linear-disc-ratio}, we have
  \( \disc(P_n^{(\lambda)};\A_{\delta},\delta) =
  \disc(P_n^{(\lambda)};S) \).
\end{exam}

\subsection{Nikolova's results}

Nikolova \cite[Section~5]{Nikolova2004} found product formulas for
discriminants and resultants involving Meixner--Pollaczek, continuous
Hahn, continuous dual Hahn, dual Hahn, Wilson, and Racah polynomials.
For Meixner--Pollaczek and continuous Hahn polynomials, the formulas
follow from \Cref{cor:linear-disc-ratio} with the shifted operator
\(Sf=(Tf)\circ v^{-1}\). For the four quadratic-lattice families, we
use \Cref{cor:avg-T-disc} and the polynomial-composition argument
given below. We first introduce the necessary definitions.

The monic \emph{Wilson polynomials} are defined by
\[
  W_n(x^2) = W_n(x^2;a,b,c,d)
  =\frac{(-1)^n(a+b)_n(a+c)_n(a+d)_n}{(n+s-1)_n}
  {}_4F_3\!\left(
    \begin{matrix}
      -n,n+s-1,a+\ii x,a-\ii x\\
      a+b,a+c,a+d
    \end{matrix};1
  \right).
\]
where \(s=a+b+c+d\). The monic \emph{Racah polynomials} are defined by
\[
  R_n(\lambda(x);\alpha,\beta,\gamma,\delta)
  =\frac{(\alpha+1)_n(\beta+\delta+1)_n(\gamma+1)_n}
          {(n+\alpha+\beta+1)_n}
  {}_4F_3\!\left(
    \begin{matrix}
      -n,n+\alpha+\beta+1,-x,x+\gamma+\delta+1\\
      \alpha+1,\beta+\delta+1,\gamma+1
    \end{matrix};1
  \right),
\]
where \(\lambda(x)=x(x+\gamma+\delta+1)\), and one of
\(\alpha+1,\beta+\delta+1,\gamma+1\) is equal to \(-N\), and
\(0\leq n\leq N\).

For a monic quadratic polynomial \( \lambda(x) \), a polynomial
\( p \) with leading coefficient \( \gamma \), and an operator
\( T \), we define
\[
  D^2(p(\lambda(x)),T) = \gamma^{-2} \res(p(\lambda(x)), T p(\lambda(x))).
\]
For a quadratic polynomial \(\lambda\) and a complex number \(h\), we
define the linear map \(\Delta_{\lambda,h}\) on polynomial
compositions \(p\circ\lambda\), with values in polynomials in \(x\), by
\[
  \bigl(\Delta_{\lambda,h}(p\circ\lambda)\bigr)(x)
  =\frac{p(\lambda(x+h))-p(\lambda(x))}
         {\lambda(x+h)-\lambda(x)}.
\]
The quotient is interpreted according to \eqref{eq:u-v}.

Using our notation, Nikolova's formulas \cite[(38),
(39)]{Nikolova2004} are written as follows:
\begin{align}
  \label{eq:nikolova-wilson}
  &D^2\bigl(W_n(x^2;a,b,c,d),\Delta_{x^2,\ii}\bigr)
  \\
  &\quad=
  \frac{1}{(n+s-1)_n^{4n-4}}
  \prod_{j=1}^n\Bigl(
    j^{2j}
    (a+b+j-1)^{2j-2}(a+c+j-1)^{2j-2}
    (a+d+j-1)^{2j-2}
  \notag\\[-2mm]
  &\qquad \times(b+c+j-1)^{2j-2}(b+d+j-1)^{2j-2}
    (c+d+j-1)^{2j-2} (s+n+j-2)^{2n-2j}
  \Bigr),
  \notag\\
  \label{eq:nikolova-racah}
  &D^2\bigl(R_n(\lambda(x);\alpha,\beta,\gamma,\delta),
    \Delta_{\lambda,1} \bigr)
  \\
  &\quad=
  \left(
    \frac{(\alpha+1)_n(\beta+\delta+1)_n(\gamma+1)_n}
         {(n+\alpha+\beta+1)_n}
  \right)^{4n-4}
  \prod_{j=1}^n\Bigl(
      j^{2j}(\beta+j)^{2j-2}
      (\alpha+\beta-\gamma+j)^{2j-2}
      (\alpha-\delta+j)^{2j-2}
  \Bigr)
  \notag\\[-1mm]
  &\qquad{}
  \times\prod_{j=1}^n
    \frac{(\alpha+\beta+n+j)^{2n-2j}}
    {(\alpha+j)^{4n-2j-2}
      (\beta+\delta+j)^{4n-2j-2}
      (\gamma+j)^{4n-2j-2}}. \notag
\end{align}
Both formulas can be obtained from \Cref{cor:avg-T-disc}. We will
briefly sketch how to obtain \eqref{eq:nikolova-wilson}.

Since resultants are invariant under translation \( x\mapsto x-i/2 \), we
rewrite the left-hand side of \eqref{eq:nikolova-wilson} as follows:
\begin{align*}
  D^2\bigl(W_n(x^2),\Delta_{x^2,\ii}\bigr)
  &=\res\left(W_n(x^2),\frac{W_n ( (x+i)^2 ) - W_n(x^2)}
    {\left(x+i\right)^2-x^2}\right)\\
  &=\res\left(W_n\left(\left(x-\frac{i}{2}\right)^2\right),
    \frac{W_n \left( \left(x+\frac{i}{2}\right)^2 \right) - W_n\left( \left(x-\frac{i}{2}\right)^2 \right)}
    {2ix}\right).
\end{align*}
By \Cref{pro:4}, adding the second component multiplied by \( ix \) to
the first component, we obtain
\[
  D^2\bigl(W_n(x^2),\Delta_{x^2,\ii}\bigr)
  = \res((\A_{\mathcal{W}} W_n)(x^2), (\mathcal{W} W_n)(x^2))
  = \res((\A_{\mathcal{W}} W_n)(x), (\mathcal{W} W_n)(x))^2,
\]
where \( \mathcal{W} \) is the Wilson operator in \eqref{eq:wilson}.
Note that the second equality above follows from \Cref{def:resultant}.
Applying \Cref{cor:avg-T-disc} yields \eqref{eq:nikolova-wilson}.

\section{Resultants and orthogonal polynomials}
\label{sec:result-orth-polyn}

In this section we give two general methods for finding two sets of
orthogonal polynomials which form a linked pair in the sense of
\Cref{def:5}. This results in two factorizations of
resultants of orthogonal polynomials (see Theorems~\ref{linearthm2}
and \ref{quadthm}). These results are independent of the orthogonal
polynomials being in the Askey and \( q \)-Askey schemes.

Recall that the Jacobi polynomials \(P_n^{(\alpha,\beta)}(x)\) and
Laguerre polynomials \(L_n^{(\alpha)}(x)\) both have classical
discriminants which factor. The derivatives
\(D(P_n^{(\alpha,\beta)}(x))\) and \(D(L_n^{(\alpha)}(x))\) are multiples
of \(P_{n-1}^{(\alpha+1,\beta+1)}(x)\) and \(L_{n-1}^{(\alpha+1)}(x),\)
respectively. The measures for the Jacobi and Laguerre polynomials are
\( (1-x)^\alpha (1+x)^\beta dx \) and \( x^\alpha e^{-x}dx \),
respectively. The measures for the derivatives
\(D(P_n^{(\alpha,\beta)}(x))\) and \(D(L_n^{(\alpha)}(x))\) correspond to
multiplying \( (1-x)^\alpha (1+x)^\beta dx \) and
\( x^\alpha e^{-x}dx \) by \((1-x)(1+x)\) and \(x\), respectively. We
shall see that multiplying the measure by a linear or quadratic
polynomial produces an OPS that forms a linked pair with the original
OPS. These measure modifications are Christoffel transformations as in
\cite[Theorem 2.7.1]{Ismail}.

Throughout this section, let
\begin{equation}
\label{3term}
p_{n+1}(x)=(x-b_n)p_n(x)-\lambda_n p_{n-1}(x),\quad p_0(x)=1, \quad p_{-1}(x)=0,
\end{equation}
be monic orthogonal polynomials with measure \( \mu(x) \).

\subsection{Shifting the measure by a linear polynomial}

We define
\begin{equation}\label{eq:1}
Q_n(x)=
\frac{\det\begin{pmatrix} 
p_n(A)&p_{n+1}(A)\\
p_n(x)&p_{n+1}(x)\\
\end{pmatrix}}{x-A}.
\end{equation}
When \(p_n(A)\ne0\) for every \(n\ge0\), the sequence
\(\{Q_n(x)\}\) is a (non-monic) OPS for the measure
\((x-A)d\mu(x)\); see \cite[Theorem 2.7.1]{Ismail}.

\begin{thm} 
\label{linearthm2}
Let \( \{p_n(x)\} = \OPS(\{b_n\},\{\lambda_n\}) \)
and let \( Q_n(x) \) be given by \eqref{eq:1}.
Then we have
\[
\res(p_{n+1}(x), Q_n(x))= (-1)^ {\binom{n+1}{2}}p_{n+1}(A)^{n}
\prod_{j=1}^n \lambda_{j}^{j}.
\]
\end{thm}

\begin{proof}
First suppose that \(p_{n+1}(A)\ne0\).
By \eqref{eq:1}, we have
\[
Q_n(x)=
\frac{p_n(A)p_{n+1}(x)-p_{n+1}(A)p_n(x)}{x-A}.
\]
If \(x_{n+1,1},\dots, x_{n+1,n+1}\) are the zeros of \(p_{n+1}(x)\),
then
\begin{align*}
\res(p_{n+1}(x), Q_n(x)) & =\prod_{j=1}^{n+1} Q_{n}(x_{n+1,j})
= \prod_{j=1}^{n+1}  \frac{-p_{n+1}(A)p_n(x_{n+1,j})}{x_{n+1,j}-A}\\
&=(-p_{n+1}(A))^{n+1}\res(p_{n+1}(x),p_{n}(x))
\prod_{j=1}^{n+1}\frac{1}{x_{n+1,j}-A}\\
&=p_{n+1}(A)^n\res(p_{n+1}(x),p_{n}(x))\\
&=(-1)^{\binom{n+1}{2}} p_{n+1}(A)^n \prod_{j=1}^n \lambda_{j}^{j},
\end{align*}
where we used the identity
\( \prod_{j=1}^{n+1}(x_{n+1,j}-A)=(-1)^{n+1}p_{n+1}(A) \) and
\Cref{thm:schur}.
This proves the formula whenever \(p_{n+1}(A)\ne0\). Since the
coefficients of \(Q_n(x)\), and hence both sides of the formula, are
polynomials in \(A\), the identity extends to every \(A\).
\end{proof}

Here is a variation of Theorem~\ref{linearthm2} with equal degree polynomials.
The proof of the following theorem is similar to that of \Cref{linearthm2}.

\begin{thm} 
\label{linearthm1}
Let \( \{p_n(x)\} = \OPS(\{b_n\},\{\lambda_n\}) \)
and let \( Q_n(x) \) be given by \eqref{eq:1}.
Then we have
\[
\res(p_{n}(x), Q_n(x))= (-1)^ {\binom{n}{2}}p_{n}(A)^{n-1}
\prod_{j=1}^n \lambda_{j}^{j}.
\]
\end{thm}

There are several examples in the Askey and \( q \)-Askey schemes when
multiplying by a linear polynomial just shifts parameters. We state
the Askey--Wilson case.

\begin{cor}
\label{2ndcor}
For \(n\geq1\), we have
\begin{align*}
  \res\bigl(p_n(x;a,b,c,d\mid q),
    p_{n-1}(x;aq,b,c,d\mid q)\bigr)
  &=(-1)^{\binom n2}
    \frac{(ab,ac,ad;q)_n^{n-1}}
         {(ab,ac,ad;q)_{n-1}^n}
    \frac{(abcdq^{n-2};q)_{n-1}^n}
         {(abcdq^{n-1};q)_n^{n-1}}
    \prod_{j=1}^{n-1}\lambda_j^{\,j},\\
  \res\bigl(p_n(x;a,b,c,d\mid q),
    p_n(x;aq,b,c,d\mid q)\bigr)
  &=(-1)^{\binom n2}
    \frac{2^na^n(abcdq^{n-1};q)_n}{(ab,ac,ad;q)_n}
    \prod_{j=1}^{n}\lambda_j^{\,j},
\end{align*}
where \( \lambda_j \) are given in \eqref{eq:AW-b,la}.
\end{cor}
\begin{proof}
  The Askey--Wilson polynomials are orthogonal with respect to the
  measure
  \[
    \mu(x;a,b,c,d\mid q)
    =\frac{(z^2,z^{-2};q)_\infty}
    {\sqrt{1-x^2}\,
    (az,az^{-1},bz,bz^{-1},cz,cz^{-1},dz,dz^{-1};q)_\infty},
  \]
  where \(x=(z+z^{-1})/2\). Note that
  \[
    \frac{\mu(x;aq,b,c,d\mid q)}{\mu(x;a,b,c,d\mid q)}
    =(1-az)(1-az^{-1})
    =1-2ax+a^2=-2a(x-A),
    \qquad A=\frac{a+a^{-1}}2.
  \]
  Thus, the polynomials \( Q_n(x) \) in \eqref{eq:1} with
  \( p_n(x) =p_n(x;a,b,c,d\mid q) \) and \( A=\frac{a+a^{-1}}2 \)
  are orthogonal with respect to the measure \( \mu(x;aq,b,c,d\mid q) \).
  Considering the leading coefficient, we obtain
  \[
    Q_n(x)=p_n(A;a,b,c,d\mid q)\,
    p_n(x;aq,b,c,d\mid q).
  \]
  By \Cref{linearthm2} and \Cref{linearthm1}, we have
  \begin{align*}
    \res\bigl(p_n(x;a,b,c,d\mid q),p_{n-1}(x;aq,b,c,d\mid q)\bigr)
    &=(-1)^{\binom n2}
    \frac{p_n(A)^{n-1}}{p_{n-1}(A)^n} \prod_{j=1}^{n-1}\lambda_j^j,\\
        \res\bigl(p_n(x;a,b,c,d\mid q),p_n(x;aq,b,c,d\mid q)\bigr)
    &=(-1)^{\binom n2}\frac{1}{p_n(A)} \prod_{j=1}^{n}\lambda_j^j.
  \end{align*}
  Since
  \[
    p_m(A;a,b,c,d\mid q)
    =\frac{a^{-m}(ab,ac,ad;q)_m}
    {2^m(abcdq^{m-1};q)_m},
  \]
  we obtain the formulas in the statement.
\end{proof}

\subsection{Shifting the measure by a quadratic polynomial}

As before, let \( p_n(x) \) be the orthogonal polynomials given by
\eqref{3term}. Let \(A\) and \(B\) be distinct constants. We define
\begin{equation}
\label{3by3det}
QQ_n(x)=
\frac{\det\begin{pmatrix} 
p_n(A)&p_{n+1}(A)&p_{n+2}(A)\\
p_n(B)&p_{n+1}(B)&p_{n+2}(B)\\
p_n(x)&p_{n+1}(x)&p_{n+2}(x)\\
\end{pmatrix}}{(x-A)(x-B)}.
\end{equation}
When
\[
p_n(A)p_{n+1}(B)-p_{n+1}(A)p_n(B)\ne0
\qquad(n\ge0),
\]
the sequence \(\{QQ_n(x)\}\) is a non-monic OPS for the measure
\((x-A)(x-B)d\mu(x)\).

\begin{thm} 
\label{quadthm}
Let \(A\ne B\), let
\( \{p_n(x)\} = \OPS(\{b_n\},\{\lambda_n\}) \),
and let \( QQ_n(x) \) be given by \eqref{3by3det}.
Then we have
\[
\res(p_{n+1}(x), QQ_n(x))= (-1)^{\binom{n+1}{2}}(B-A)^{n+1}p_{n+1}(A)^n p_{n+1}(B)^n
\prod_{j=1}^n \lambda_j^{j}.
\]
\end{thm}

\begin{proof}
First suppose that \(p_{n+1}(A)p_{n+1}(B)\ne0\).
From \eqref{3by3det} we have
\[
(x-A)(x-B)QQ_{n}(x)=D_n p_{n+2}(x)+E_n p_{n+1}(x)+F_n p_{n}(x),
\]
where \( E_n \) is a certain constant and
\begin{align*}
D_n&=p_n(A)p_{n+1}(B)-p_{n+1}(A)p_n(B), \\
F_n&=p_{n+1}(A)p_{n+2}(B)-p_{n+2}(A)p_{n+1}(B).
\end{align*}
Let \(x_{n+1,1},\dots, x_{n+1,n+1}\) be the zeros of
\(p_{n+1}(x)\). Then
\begin{align*}
\res(p_{n+1}(x), QQ_n(x))&= \prod_{j=1}^{n+1} QQ_{n}(x_{n+1,j})
= \prod_{j=1}^{n+1} \frac{D_n p_{n+2}(x_{n+1,j})+F_n p_{n}(x_{n+1,j})}{(x_{n+1,j}-A)(x_{n+1,j}-B)}\\
&= \frac{1}{p_{n+1}(A)p_{n+1}(B)}
\prod_{j=1}^{n+1}\left(
D_n p_{n+2}(x_{n+1,j})+F_n p_{n}(x_{n+1,j})\right).
\end{align*}
The three-term recurrence \eqref{3term}, with \(n\) replaced by
\(n+1\), implies
\(p_{n+2}(x_{n+1,j})=-\lambda_{n+1} p_{n}(x_{n+1,j}) \). Thus
\begin{align*}
\res(p_{n+1}(x), QQ_n(x))&=\frac{1}{p_{n+1}(A)p_{n+1}(B)} 
\prod_{j=1}^{n+1}\left(
(-\lambda_{n+1}D_n+F_n)p_{n}(x_{n+1,j})\right)\\
&=\frac{(-\lambda_{n+1}D_n+F_n)^{n+1}}{p_{n+1}(A)p_{n+1}(B)}
\res(p_{n+1}(x), p_n(x))\\
&=\frac{(-\lambda_{n+1}D_n+F_n)^{n+1}}{p_{n+1}(A)p_{n+1}(B)}
(-1)^{\binom{n+1}{2}}\prod_{j=1}^n \lambda_j^j \qquad (\mbox{by \Cref{thm:schur}}).
\end{align*}
Since
\begin{align*}
-\lambda_{n+1}D_n+F_n
&=-\lambda_{n+1}\left(p_n(A)p_{n+1}(B)-p_{n+1}(A)p_n(B)\right)
+p_{n+1}(A)p_{n+2}(B)-p_{n+2}(A)p_{n+1}(B)\\
&=p_{n+1}(B)(-\lambda_{n+1}p_n(A)-p_{n+2}(A))+p_{n+1}(A)(\lambda_{n+1}p_n(B)+p_{n+2}(B))\\
&=p_{n+1}(B)(-(A-b_{n+1})p_{n+1}(A))+p_{n+1}(A)((B-b_{n+1})p_{n+1}(B)) \\
&=(B-A) p_{n+1}(A)p_{n+1}(B),
\end{align*}
we obtain the identity whenever
\(p_{n+1}(A)p_{n+1}(B)\ne0\). The determinant in
\eqref{3by3det} is divisible by \((x-A)(x-B)\), so the coefficients of
\(QQ_n(x)\), and hence both sides of the formula, are polynomials in
\(A\) and \(B\). Therefore the identity extends to all \(A\ne B\).
\end{proof}

Using \Cref{quadthm}, we obtain an explicit product formula for the
following resultant involving Askey--Wilson polynomials.

\begin{cor}
\label{cor:aw-quadratic-shift}
For \(n\geq1\), we have
\begin{align*}
  &\res\bigl(p_n(x;a,b,c,d\mid q),
    p_{n-1}(x;aq,bq,c,d\mid q)\bigr)\\
  &=(-1)^{\binom n2}
    \frac{
      (1-ab)^n(1-abcdq^{2n-3})^n
      (ab,ac,ad;q)_n^{n-1}(ab,bc,bd;q)_n^{n-1}
      (abcdq^{n-2};q)_{n-1}^{2n}
    }{
      (1-abq^{n-1})^n(1-abcdq^{n-2})^n
      (ab,ac,ad;q)_{n-1}^n(ab,bc,bd;q)_{n-1}^n
      (abcdq^{n-1};q)_n^{2n-2}
    }
    \prod_{j=1}^{n-1}\lambda_j^{\,j},
\end{align*}
where \(\lambda_j\) are given in \eqref{eq:AW-b,la}.
\end{cor}
\begin{proof}
  Let \( A=(a+a^{-1})/2 \) and \( B=(b+b^{-1})/2 \). The
  Askey--Wilson measures satisfy
  \begin{align*}
    \frac{\mu(x;aq,bq,c,d\mid q)}{\mu(x;a,b,c,d\mid q)}
    &=(1-az)(1-az^{-1})(1-bz)(1-bz^{-1})\\
    &=(1-2ax+a^2)(1-2bx+b^2)
      =4ab(x-A)(x-B).
  \end{align*}
  Thus, the polynomials \(QQ_n(x)\) in \Cref{3by3det} with
  \( p_n(x) = p_n(x;a,b,c,d\mid q) \) are orthogonal with respect to
  the measure \( \mu(x;aq,bq,c,d\mid q) \).
  Considering the leading coefficient, we have
    \[
      QQ_m(x)=\Delta_m p_m(x;aq,bq,c,d\mid q),
  \]
  where \( \Delta_m =p_m(A)p_{m+1}(B)-p_{m+1}(A)p_m(B) \).
  Hence, by the same argument as in the proof of \Cref{2ndcor},
  we obtain the desired formula.
\end{proof}

Koornwinder and Bouzeffour \cite[(4.14)]{Koornwinder2011} studied the
operator \(Y_{21}\), which sends \(p_n(x;a,b,c,d|q)\) to a scalar
multiple of \(p_{n-1}(x;aq,bq,c,d|q)\). The operator \(Y_{21}\) can be
written as a linear combination of the identity operator and the
Askey--Wilson averaging and difference operators with base
\( q^2 \), whose coefficients are rational functions in \(x\). We do not give the
details.

We give applications of \Cref{quadthm} for \(q\)-Racah and Pollaczek
polynomials. We note that the Pollaczek polynomials are not in the
Askey scheme. \Cref{quadthm} may be applied to any of the
polynomials in the Askey or \(q\)-Askey scheme with an appropriate
choice of \(A\) and \(B\).

The monic \emph{\(q\)-Racah polynomials} are defined by
\[
  R_n(\mu(x;\gamma,\delta);\alpha,\beta,\gamma,\delta\mid q)
  =\frac{(\alpha q,\beta\delta q,\gamma q;q)_n}
    {(\alpha\beta q^{n+1};q)_n}
  {}_4\phi_3\!\left(
    \begin{matrix}
      q^{-n},\alpha\beta q^{n+1},q^{-x},\gamma\delta q^{x+1}\\
      \alpha q,\beta\delta q,\gamma q
    \end{matrix};q,q
  \right) \qquad (0\le n\le N),
\]
where
\( X=\mu(x)=\mu(x;\gamma,\delta)=q^{-x}+\gamma\delta q^{x+1} \).
Here, \( N \) is a fixed integer, and one of \(\alpha q\),
\(\beta\delta q\), and \(\gamma q\) is equal to \(q^{-N}\). The
\( q \)-Racah polynomials are orthogonal with respect to the discrete
measure
\[
  w(x;\alpha,\beta,\gamma,\delta)
  =
  \frac{(\alpha q,\beta\delta q,\gamma q,\gamma\delta q;q)_x}
       {(q,\alpha^{-1}\gamma\delta q,\beta^{-1}\gamma q,\delta q;q)_x}
       \frac{1-\gamma\delta q^{2x+1}}{(\alpha\beta q)^x(1-\gamma\delta q)},
\]
with support \( x=0,1,\dots,N \).
We consider the case \(\alpha=q^{-1-N}\).

\begin{cor}
For \(1\leq n\leq N-1\), we have
\begin{align*}
  &\res\bigl(
    R_n(\mu(x;\gamma,\delta);q^{-1-N},\beta,\gamma,\delta\mid q),
    R_{n-1}(\mu(x-1;\gamma q,\delta q);
      q^{1-N},\beta,\gamma q,\delta q\mid q)\bigr)\\
  &=
    \res\bigl(
    R_n(\mu(x;\gamma,\delta);q^{-1-N},\beta,\gamma,\delta\mid q),
    R_{n-1}(q\mu(x;\gamma,\delta);
      q^{1-N},\beta,\gamma q,\delta q\mid q)\bigr)\\
  &=(-1)^{\binom n2}
    q^{n(n-1)}
    \left(
    \frac{\gamma q^{n-N}(1-q^N)
      (1-\delta q^{N-n+1})
      (1-\beta\gamma^{-1}q^{n-1-N})}
    {1-\beta q^{2n-N-1}}
    \right)^n
    \left(\frac{\kappa_{n-1}}{\kappa_n}\right)^n
    \frac{\kappa_n^2}{E_n}
    \prod_{j=1}^{n-1}\lambda_j^{\,j},
\end{align*}
where both resultants are taken with respect to
\(X=\mu(x;\gamma,\delta)\), and
\begin{align*}
  \kappa_m
  &=\frac{(\beta q^{m-N};q)_m}
    {(q^{-N},\beta\delta q,\gamma q;q)_m},
    \qquad E_m =\frac{(\delta q^{N-m+1},\beta\gamma^{-1}q^{-N};q)_m}
    {(\beta\delta q,\gamma^{-1}q^{-m};q)_m},
    \qquad \lambda_j=\mathsf A_{j-1}\mathsf C_j,\\ 
  \mathsf A_j
  &=
  \frac{(1-q^{j-N})(1-\beta q^{j-N})
    (1-\beta\delta q^{j+1})(1-\gamma q^{j+1})}
  {(1-\beta q^{2j-N})(1-\beta q^{2j+1-N})},\\
  \mathsf C_j
  &=
  \frac{q(1-q^j)(1-\beta q^j)
    (\gamma-\beta q^{j-1-N})(\delta-q^{j-1-N})}
  {(1-\beta q^{2j-1-N})(1-\beta q^{2j-N})}.
\end{align*}
\end{cor}

\begin{proof}
  Since \( \mu(x-1;\gamma q,\delta q)=q\mu(x;\gamma,\delta) \), the
  two resultants are equal. Let
  \(A=\mu(0;\gamma,\delta)=1+\gamma\delta q\), and
  \(B=\mu(N;\gamma,\delta)=q^{-N}+\gamma\delta q^{N+1}\). For
  \(1\leq x\leq N-1\), we have
  \begin{align}
    \label{eq:6}    \frac{w(x-1;q^{1-N},\beta,\gamma q,\delta q)}
    {w(x;q^{-1-N},\beta,\gamma,\delta)}
    &=c\, (1-q^{-x})(1-\gamma\delta q^{x+1})
      (1-q^{N-x})(1-\gamma\delta q^{N+x+1})\\
    \notag &=cq^N(X-A)(X-B),
  \end{align}
  where
  \[
    c=
    \frac{\beta(1-\beta^{-1}\gamma q)(1-\delta q)}
    {q^{2N-2}(1-q^{-N})(1-q^{1-N})
      (1-\beta\delta q)(1-\gamma q)
      (1-\gamma\delta q^2)(1-\gamma\delta q^3)}.
  \]
  Thus the polynomials \(QQ_m(X)\) in \eqref{3by3det} with
  \( p_n(x) = R_n(\mu(x;\gamma,\delta);\alpha,\beta,\gamma,\delta\mid
  q) \) are orthogonal with respect to the measure
  \(w(x-1;q^{1-N},\beta,\gamma q,\delta q)\) with support
  \( x=0,1,\dots,N \). By \eqref{eq:6}, this measure vanishes at
  \( x=0 \) and \( x=N \), hence, we may assume that the support is
  \( x=1,\dots,N-1 \).

  On the other hand, since
  \( \mu(x-1;\gamma q,\delta q)=q\mu(x;\gamma,\delta) = qX \), the polynomials
  \( R_m(qX;q^{1-N},\beta,\gamma q,\delta q\mid q) \) are orthogonal
  with respect to the same measure
  \(w(x-1;q^{1-N},\beta,\gamma q,\delta q)\) with support
  \( x=1,\dots,N-1 \). Hence, considering the leading coefficient,
  we obtain
  \[
    QQ_m(X)=\frac{\Delta_m}{q^m}
    R_m(qX;q^{1-N},\beta,\gamma q,\delta q\mid q),
    \qquad
    \Delta_m=R_m(A)R_{m+1}(B)-R_{m+1}(A)R_m(B).
  \]
  The terminating \(q\)-Pfaff--Saalsch\"utz summation gives
  \begin{align*}
    R_m(A)
    &=\frac{(q^{-N},\beta\delta q,\gamma q;q)_m}
      {(\beta q^{m-N};q)_m},\\
    R_m(B)
    &=\frac{(q^{-N},\gamma q,\delta q^{N-m+1},
      \beta\gamma^{-1}q^{-N};q)_m}
      {(\beta q^{m-N},\gamma^{-1}q^{-m};q)_m}.
  \end{align*}
  By \Cref{quadthm}, we obtain the desired formula.
\end{proof}

The monic \emph{Pollaczek polynomials}
\(P_n^{(\lambda)}(x;a,b)\) (see \cite[\S 18.35]{NIST:DLMF}) are given by
\begin{equation}\label{eq:7}
P_n^{(\lambda)}(x;a,b) =
\frac{n!}{2^n(a+\lambda)_n}
\sum_{k=0}^n
 \frac{(\lambda+\ii\tau_{a,b}(\theta))_k}{k!}
 \frac{(\lambda-\ii\tau_{a,b}(\theta))_{n-k}}{(n-k)!}
e^{\ii(n-2k)\theta},
\end{equation}
where \(x=\cos\theta\), \((z)_k=z(z+1)\cdots(z+k-1)\), and
\(\tau_{a,b}(\theta)=\frac{ax+b}{\sin\theta} \). They are
hypergeometric polynomials but are not in the Askey scheme. They
satisfy the three-term recurrence relation
\[
  P_{n+1}^{(\lambda)}(x;a,b)
  =\left(x+\frac{b}{n+a+\lambda}\right)P_n^{(\lambda)}(x;a,b)
  -\frac{n(n+2\lambda-1)}
  {4(n+a+\lambda-1)(n+a+\lambda)}
  P_{n-1}^{(\lambda)}(x;a,b),
\]
with \( P_{-1}^{(\lambda)}(x;a,b)=0 \) and
\( P_0^{(\lambda)}(x;a,b)=1 \). For \(\lambda>0\) and
\(a\geq b\geq-a\), they are orthogonal on \((-1,1)\) with respect to
the absolutely continuous measure \(w^{(\lambda)}(x;a,b)\,dx\), where,
for \(x=\cos\theta\) and \(0<\theta<\pi\),
\[
  w^{(\lambda)}(\cos\theta;a,b)
  =\frac{1}{\pi}e^{(2\theta-\pi)\tau_{a,b}(\theta)}
  (2\sin\theta)^{2\lambda-1}
  \left|\Gamma\bigl(\lambda+\ii\tau_{a,b}(\theta)\bigr)\right|^2.
\]

\begin{cor}
  We have
\[
\res(P_{n+1}^{(\lambda)}(x;a,b), P_{n}^{(\lambda+1)}(x;a,b)) =
(-1)^{\binom{n+1}{2}}
\left(\frac{\lambda}{a+\lambda}\right)^n
\frac{n!}{2^n(a+\lambda+1)_n}
\prod_{k=1}^{n-1}
\left(
\frac{(k+1)(2\lambda+k+1)}
{4(a+\lambda+k+1)^2}
\right)^k.
\]
\end{cor}
\begin{proof}
For generic parameters satisfying \(\lambda>0\) and \(a>b>-a\), we
have
\begin{align}\label{eq:9}
  \frac{w^{(\lambda+1)}(\cos\theta;a,b)}
       {w^{(\lambda)}(\cos\theta;a,b)}
  &=(2\sin\theta)^2 \left|\lambda+\ii\tau_{a,b}(\theta)\right|^2
  =(2\sin\theta)^2 (\lambda+\ii\tau_{a,b}(\theta))(\lambda-\ii\tau_{a,b}(\theta))\\
  \notag
    &=4\left(\lambda^2(1-x^2)+(ax+b)^2\right)
    = 4(a^2-\lambda^2)(x-A)(x-B),
\end{align}
for some constants \(A\) and \(B\). Let \(QQ_m(x)\) be the polynomial
in \eqref{3by3det} with \( p_m(x)=P_m^{(\lambda)}(x;a,b) \).
Then 
\[
  QQ_m(x)=\Delta_m P_m^{(\lambda+1)}(x;a,b),
  \qquad
  \Delta_m=p_m(A)p_{m+1}(B)-p_{m+1}(A)p_m(B).
\]
Hence, by \Cref{quadthm}, the resultant in the statement is equal to
\[
  \Delta_n^{-(n+1)}\res(p_{n+1}(x),QQ_n(x))
  =(-1)^{\binom{n+1}{2}}
    \frac{(B-A)^{n+1}p_{n+1}(A)^np_{n+1}(B)^n}
         {\Delta_n^{n+1}} \prod_{j=1}^n\Lambda_j^j,
\]
where
\[
  \Lambda_m=
  \frac{m(m+2\lambda-1)}
       {4(m+a+\lambda-1)(m+a+\lambda)}.
\]

It remains to evaluate the factors involving \(A\) and \(B\). By
\eqref{eq:9}, we can take \( \theta_1 \) and \( \theta_2 \) so that
\( A = \cos\theta_1 \), \( B = \cos\theta_2 \), and
\(
(\lambda+\ii\tau_{a,b}(\theta_r))(\lambda-\ii\tau_{a,b}(\theta_r))=0
\) for \( r=1,2 \). Since
\( \tau_{a,b}(-\theta_r) = -\tau_{a,b}(\theta_r) \), replacing
\( \theta_r \) by \( -\theta_r \) if necessary, we may assume that
\( \lambda+\ii \tau_{a,b}(\theta_r)=0 \). Let
\( z_r = e^{\ii\theta_r} \). Then, by \eqref{eq:7}, we obtain
\[
  p_m(A)=c_mz_1^m,
  \qquad
  p_m(B)=c_mz_2^m,
  \qquad 
  c_m=\frac{(2\lambda)_m}{2^m(a+\lambda)_m}.
\]
For \(z=e^{\ii\theta}\), the condition \( \lambda+\ii \tau_{a,b}(\theta)=0 \) is
equivalent to
\[
  (a-\lambda)z^2+2bz+a+\lambda=0.
\]
Since \(z_1,z_2\) are the roots of this equation, we have
\[
  z_1z_2=\frac{a+\lambda}{a-\lambda}, \qquad 
  B-A =\frac{\lambda}{a+\lambda}(z_2-z_1), \qquad 
  \Delta_m =c_mc_{m+1}
  \left(\frac{a+\lambda}{a-\lambda}\right)^m(z_2-z_1).
\]
Combining these results, we obtain the desired formula.
\end{proof}

It is possible to prove Theorems~\ref{linearthm2} and \ref{quadthm}
using \Cref{lem:key}. To this end, we need the following
proposition.

\begin{prop}\label{pro:2}
  Let \(A\ne B\), let
  \( \{p_n(x)\} = \OPS(\{b_n\},\{\lambda_n\}) \), and let
  \( Q_n(x) \) and \( QQ_n(x) \) be given by \eqref{eq:1} and
  \eqref{3by3det}. Assume, for every \(n\ge0\), that
  \[
    p_n(A)\ne0
    \quad\text{and}\quad
    p_n(A)p_{n+1}(B)-p_{n+1}(A)p_n(B)\ne0.
  \]
  Then there exists constants \(A_n, B_n,C_n, E_n\)
  and \(F_n\) such that
\begin{align*}
p_n(x)&=A_nQ_n(x)+B_nQ_{n-1}(x),\\
p_n(x)&= C_nQQ_n(x)+E_nQQ_{n-1}(x)+F_nQQ_{n-2}(x).
\end{align*}
Consequently, \( (\{p_n(x)\}, \{\widehat Q_n(x)\}) \) and
\( (\{p_n(x)\}, \{\widehat{QQ}_n(x)\}) \) are linked pairs, where
\( \widehat Q_n(x) \) and \( \widehat{QQ}_n(x) \) denote \( Q_n(x) \)
and \( QQ_n(x) \) divided by their leading coefficients: substituting
the three-term recurrence of \( \widehat Q_n(x) \)
(resp.\ \( \widehat{QQ}_n(x) \)) into the above identities puts them
in the form required by \Cref{def:5}.
\end{prop}

\begin{proof}
  We prove the statement for \(Q_n(x)\), the \(QQ_n(x)\) proof is
  nearly identical. Let
\begin{equation}\label{eq:2}
p_n(x)=\sum_{k=0}^n c_{n,k} Q_k(x).
\end{equation}
For each \(k\),
\[
Q_k(x)=
\frac{p_k(A)p_{k+1}(x)-p_{k+1}(A)p_k(x)}{x-A}.
\]
The numerator is a polynomial in \(A\) and \(x\), and it vanishes when
\(x=A\). Hence it is divisible by \(x-A\), so the coefficients of
\(Q_k(x)\), as a polynomial in \(x\), are polynomials in \(A\).
Therefore, in the expansion \eqref{eq:2}, each coefficient \(c_{n,k}\)
is a rational function of \(A\).

Let \( k<n-1 \). It suffices to show that \(c_{n,k}=0\). By the
orthogonality of \(Q_n\) we have
\[
 c_{n,k} \int_{-\infty}^\infty Q_k(x) Q_k(x) (x-A) d\mu(x) =
 \int_{-\infty}^\infty p_n(x) Q_k(x) (x-A) d\mu(x)=0,
\]
where the last equality follows from the orthogonality of \(p_n(x)\)
since \( k<n-1 \). The integral multiplying \(c_{n,k}\) is zero only
for a finitely many values of \(A\). Thus, \(c_{n,k}=0\) for all but
a finite many values of \(A\). Since \(c_{n,k}\) is a rational
function of \(A\), it is identically zero, which completes the proof.
\end{proof}

Once \( QQ_n(x) \) is expressed in terms of \( QQ_{n-1}(x) \) and
\( QQ_{n-2}(x) \), the identities in \Cref{pro:2} give mixed
three-term relations of the form in \Cref{lem:key}. One can find
\(A_n, B_n, C_n, E_n\), and \(F_n\) explicitly, apply
\Cref{lem:key}, and thereby obtain Theorems~\ref{linearthm2}
and~\ref{quadthm}. We omit the details.

\subsection{More variations}

For completeness we state, without proof, the variations of the key
theorems for type \(R_I\) and type \(R_{II}\) orthogonal polynomials.
For our purpose, one may consider type \( R_{II} \) orthogonal
polynomials as polynomials satisfying the recurrence
\[
  P_{n+1}(x)= (x-b_n)P_n(x)-(c_nx^2+ a_n x+\lambda_n)P_{n-1}(x), \quad
  P_0(x)=1, \quad P_{-1}(x)=0.
\]
If \( c_n=0 \), then they are called type \( R_I \) orthogonal
polynomials. We refer the reader to \cite{IsmailMasson} for more
details.

In the resultant formulas below, for
\(P_k(x)=\sum_{j=0}^k p_{k,j}x^j\), we interpret
\[
  A^kP_k(-\Lambda/A)
  :=\sum_{j=0}^k p_{k,j}(-\Lambda)^j A^{k-j}.
\]
This polynomial interpretation also covers \(A=0\), when the
expression equals \(p_{k,k}(-\Lambda)^k\).

Here is the type \(R_I\) version of Schur's lemma, \Cref{thm:schur}.

\begin{prop} Let \(P_n(x)\) be monic type \(R_I\) orthogonal polynomials,
\label{typeRIthm}
\begin{equation}
\label{typeRI}
P_{n+1}(x)= (x-B_n)P_n(x)-(A_n x+\Lambda_n)P_{n-1}(x), \quad P_0(x)=1, \quad P_{-1}(x)=0.
\end{equation}
Then 
\[
\res(P_{n+1}(x),P_{n}(x))= \prod_{k=1}^n  A_k^k P_k(-\Lambda_k/A_k).
\]
\end{prop}

Here is the type \(R_I\) version of \Cref{lem:key}.

\begin{prop} Suppose that \(P_n(x)\) is a monic type \(R_I\) OPS
  satisfying \eqref{typeRI}, and for a fixed \(n\ge 2\) and constants
  \(\beta_n\) and \(\Gamma_n\), let
\[
G_n(x)= (x-\beta_n)P_{n-1}(x) -\Gamma_nP_{n-2}(x).
\]
Then
\[
\res(G_n(x),P_{n-1}(x)) = (-\Gamma_n)^{n-1} \prod_{k=1}^{n-2}  A_k^k\  P_k(-\Lambda_k/A_k).
\]
\end{prop}

Here is the type \(R_{II}\) version of Schur's lemma, \Cref{thm:schur}.

\begin{prop} 
\label{SchurRII}
Suppose that \(P_n(x)\), \( n\ge 0 \), is a polynomial of degree \(n\) satisfying
\begin{equation}
\label{genrecur}
P_{n+1}(x)= L_n(x)P_n(x)- (A_{n,1} x+\Lambda_{n,1})(A_{n,2} x+\Lambda_{n,2})P_{n-1}(x), \quad P_0(x)=1, \quad P_{-1}(x)=0
\end{equation}
for some linear polynomials \(L_n(x)\). Then 
\[
\res(P_{n+1}(x),P_{n}(x))= (-1)^{\binom{n+1}{2}}\prod_{k=1}^n\prod_{t=1}^2  A_{k,t}^k P_k(-\Lambda_{k,t}/A_{k,t}).
\]
\end{prop}

Here is the type \(R_{II}\) version of \Cref{lem:key}. The leading
coefficient of a polynomial \( f \) is denoted by \( \lc (f) \).

\begin{prop}\label{pro:3} 
Suppose that \(P_n(x)\) is a set of polynomials satisfying \eqref{genrecur}, and 
for a fixed \(n\ge 2\), let
\[
G_n(x)= (x-B_n)P_{n-1}(x) -\Gamma_nP_{n-2}(x).
\]
Then
\[
  \res(G_n(x),P_{n-1}(x)) = (-1)^{\binom{n-1}{2}}(-\Gamma_n)^{n-1}
  \lc(P_{n-1})^2\prod_{k=1}^{n-2}\prod_{t=1}^2 A_{k,t}^k
  P_k(-\Lambda_{k,t}/A_{k,t}).
\]
\end{prop}

There are versions of these theorems involving higher-degree
polynomials, but we do not state them.

\section{Companion OPS}
\label{sec:companion-ops}

In this section, we introduce a finite OPS with one polynomial
appended, called an \( (f,g) \)-companion OPS (or \( (f,g) \)-COPS)
for given polynomials \( f \) and \( g \). We show that the resultant
\( \res(f,g) \) can be expressed as a scalar multiple of the Hankel
determinant of the moments of the \( (f,g) \)-COPS. We give explicit constructions of
\((f_N,Tf_N)\)-COPS and \((\A_{u,v}f_N,\D_{u,v}f_N)\)-COPS under
certain conditions, which are applicable to all orthogonal polynomials
in the Askey and \( q \)-Askey schemes.

We briefly explain our motivation to consider COPS. Suppose that we
want to compute the discriminant of
\( f(x) = (x-\alpha_1) \cdots (x-\alpha_N) \). Write the power-sum
symmetric polynomial
\( p_k = p_k(\alpha_1, \dots, \alpha_N) = \sum_{ j=1 }^{ N }
\alpha_j^k \) of the zeros of \( f \), for \( k \geq 1 \), with
\( p_0 = N \), and let \( V = (\alpha_i^{j-1})_{i,j = 1}^N \). Then we
have
\begin{equation}\label{eq:12}
  \disc(f) = \prod_{1\le i<j\le N} (\alpha_i-\alpha_j)^2 =
  \det(V)^2= \det(V^T V) = \det(p_{i+j})_{i,j=0}^{N-1}.
\end{equation}
Suppose that \( \{Q_n(x)\}= \OPS(\{b_n\},\{\lambda_n\}) \) is an OPS
whose moments satisfy \( \mu_k = p_k/N \) for \( 0\le k\le 2N-2 \).
Then using \eqref{eq:12} and the well-known fact on the Hankel
determinant of moments, we can evaluate the discriminant:
\begin{equation}\label{eq:17}
  \disc(f) = N^N \det \left( \mu_{i+j} \right)_{i,j = 0}^{N-1}
  = N^N \lambda_1^{N-1} \lambda_2^{N-2} \cdots \lambda_{N-1}^1.
\end{equation}
This motivates us to find such an OPS, which can be achieved by
finding the \( (f,f') \)-COPS. Similarly, the resultant
\( \res(f,g) \) can be computed using the \( (f,g) \)-COPS.

Before defining COPS, we introduce extended OPS and their properties.

\subsection{Extended OPS}

Recall that \(\lc(h)\) denotes the leading coefficient of a polynomial
\(h\). For a nonzero polynomial \( h \) of degree \( d \), its
\emph{reciprocal polynomial} \( h^*(x) \) is defined by
\[
  h^*(x)=x^d h(x^{-1}).
\]
Suppose \(\{Q_n(x):0\leq n\leq N\}\) is a sequence of polynomials
satisfying
\[
  Q_{n+1}(x)=(x-b_n)Q_n(x)-\lambda_nQ_{n-1}(x),
  \qquad 0\leq n<N,
\]
where \( Q_{-1}(x)=0 \) and \(Q_0(x)=1\). We define the \emph{first
  associated polynomials} \( \{Q_{n}^{(1)}(x):0\leq n\leq N-1\} \) by
\[
  Q_{n+1}^{(1)}(x)
  =(x-b_{n+1})Q_n^{(1)}(x)-\lambda_{n+1}Q_{n-1}^{(1)}(x),
  \qquad 0\leq n<N-1,
\]
where \( Q^{(1)}_{-1}(x)=0 \) and \(Q^{(1)}_0(x)=1\).

The following definition will be useful to state our results.

\begin{defn}
  Let \( \LL \) be a linear functional on the space of polynomials and
  let \( N \) be a positive integer. A sequence
  \( \{Q_n(x):0\le n\le N \} \) of polynomials with \(\deg Q_n=n\) is
  called an \emph{extended orthogonal polynomial sequence (extended
    OPS)} with respect to \( \LL \) if the following conditions hold:
  \begin{enumerate}
  \item \( \LL(Q_n(x)Q_m(x))=0 \) for \( 0\le n\ne m <N \), and
  \( \LL(Q_n(x)^2)\ne 0 \) for \( 0\le n <N \).
  \item \( \LL(x^k Q_N(x))=0 \) for all \( k\ge0 \).
  \end{enumerate}
\end{defn}

Observe that an extended OPS \( \{f_n(x):0\le n\le N \} \) is a finite
OPS \( \{f_n(x):0\le n< N \} \) with the polynomial \(f_N\) appended.
The concept of extended OPS has been studied by Maeda and Tsujimoto
\cite[Section~2.2, (2.13a)--(2.13b)]{MaedaTsujimoto2016}.

\begin{prop}\label{prop:finite-moment-next}
  Let \(N\geq1\), and let \(\mu_0,\dots,\mu_{2N-1}\) be complex
  numbers with \( \mu_0=1 \). Suppose that
  \[
    \Delta_k:=\det(\mu_{i+j})_{i,j=0}^{k-1}\neq0,
    \qquad 0\leq k\leq N,
  \]
  where \( \Delta_0=1 \). Then the following hold.
  \begin{enumerate}
  \item There is a unique linear functional \(\LL\) on the space of
    all polynomials and a unique monic extended OPS
    \( \{Q_n(x):0\le n\le N \} \) with respect to \(\LL\) such that
  \[
    \LL(x^n)=\mu_n\quad(0\leq n\leq2N-1).
  \]
\item The extended OPS satisfies
  \begin{equation}\label{eq:19}
    Q_{n+1}(x)=(x-b_n)Q_n(x)-\lambda_nQ_{n-1}(x),
    \qquad 0\leq n\le N-1,
  \end{equation}
  where \(Q_{-1}(x)=0\), \( Q_0(x)=1 \),
  \(h_n:=\LL(Q_n^2)=\Delta_{n+1}/\Delta_n\) for \( 0\leq n\le N-1 \),
  and
  \[
    b_n=\frac{\LL(xQ_n^2)}{h_n}\quad(0\leq n\le N-1),
    \qquad
    \lambda_n=\frac{h_n}{h_{n-1}}
    =\frac{\Delta_{n+1}\Delta_{n-1}}{\Delta_n^2}
    \quad(1\leq n\le N-1).
  \]
  In particular, \( \lambda_n\ne 0 \) for \( 1\le n\le N-1 \) and
  \begin{equation}\label{eq:18}
  \Delta_k=\prod_{j=1}^{k-1}\lambda_j^{k-j},
  \qquad 1\leq k\leq N.
\end{equation}

\item The moment generating function is given by
  \[
    \sum_{n\geq0}\LL(x^n)z^n
    =\frac{(Q_{N-1}^{(1)})^*(z)}{Q_N^*(z)}.
  \]
  In particular, if \(Q_N\) has \(N\) distinct zeros
  \(\alpha_1,\dots,\alpha_N\), then \(\LL\) is represented by the
  discrete measure
  \[
    \LL(x^n)=\sum_{i=1}^N w_i\alpha_i^n,
    \qquad
    w_i=\frac{Q_{N-1}^{(1)}(\alpha_i)}{Q_N'(\alpha_i)},
    \qquad n\geq0.
  \]
  \end{enumerate}
\end{prop}

\begin{proof}
  Let \(\LL_0\) be the linear functional on polynomials of degree at
  most \(2N-1\) determined by \(\LL_0(x^j)=\mu_j\). By the same
  arguments in the proofs of \cite[Chapter~I, Theorems~3.1, 3.2, 4.1,
  and 4.2]{Chihara}, we obtain that there is a unique monic OPS
  \( Q_0,\dots,Q_{N-1} \) with respect to \( \LL_0 \), and it
  satisfies the statements in Item (2), except the recurrence
  \eqref{eq:19} at \( n=N-1 \).

  Now suppose that \( Q_N=x^N + \sum_{j=0}^{N-1} c_j Q_j \) is a monic
  polynomial of degree \( N \) satisfying \( \LL_0(x^j Q_N) =0 \) for
  all \( 0\le j\le N-1 \). Since \( \LL_0(Q_j Q_N)=0 \) for all
  \( 0\le j\le N-1 \), by the orthogonality, we obtain
  \( c_j = -\LL_0(x^NQ_j)/\LL_0(Q_j^2) \). Thus such \( Q_N \) exists
  and it is unique. There is a unique extension \( \LL \) of
  \( \LL_0 \) to all polynomials satisfying \( \LL(x^k Q_N) = 0 \) for
  \( k\ge N \). This shows Item (1).
  The recurrence \eqref{eq:19} at \( n=N-1 \) can be proved similarly
  as in the previous case \( 0\le n\le N-2 \).

  For Item (3), let \(J = (J_{i,j})_{i,j=0}^{N-1}\) be the
  \( N\times N \) tridiagonal matrix with diagonal entries
  \(b_0,\dots,b_{N-1}\), subdiagonal entries all \(1\), and
  superdiagonal entries \(\lambda_1,\dots,\lambda_{N-1}\). By
  \eqref{eq:19}, the matrix \( J \) represents multiplication by \(x\)
  modulo \(Q_N\) in the basis \(Q_0,\dots,Q_{N-1}\). Note that if
  \( x^n \equiv \sum_{j=0}^{N-1} d_j Q_j \pmod{Q_N} \), then
  \( \LL(x^n) = d_0=(J^n)_{0,0}\). Thus
  \[
    \sum_{n\geq0}\LL(x^n)z^n
    = \sum_{n\geq0} (J^n)_{0,0}z^n
    = ((I-zJ)^{-1})_{0,0},
  \]
  where \( I \) is the \( N\times N \) identity matrix. On the other
  hand, by \eqref{eq:19}, we have \( Q_N(z) = \det(zI-J) \), or
  equivalently, \( Q^*_N(z) = \det(I-zJ) \). Similarly,
  \( (Q^{(1)}_{N-1})^*(z) \) is the determinant of the matrix
  \( I-zJ \) with the \( 0 \)th row and column removed. Therefore, the
  \( (0,0) \)-entry of the inverse of \( (I-zJ) \) is equal to
  \[
    ((I-zJ)^{-1})_{0,0} =\frac{(Q_{N-1}^{(1)})^*(z)}{Q_N^*(z)}.
  \]
  This shows the first statement of Item (3). If \(Q_N\) has distinct
  zeros \(\alpha_1,\dots,\alpha_N\), partial fractions give
  \[
    \frac{(Q_{N-1}^{(1)})^*(z)}{Q_N^*(z)}
    =\sum_{i=1}^N
      \frac{Q_{N-1}^{(1)}(\alpha_i)}{Q_N'(\alpha_i)}
      \frac{1}{1-\alpha_i z},
  \]
  which shows the last statement.
\end{proof}

The converse of \Cref{prop:finite-moment-next} also holds.

\begin{prop}\label{prop:finite-moment-converse}
  Let \(N\geq1\), and let \(Q_0,\dots,Q_N\) be the monic polynomials
  defined by
  \[
    Q_{n+1}(x)=(x-b_n)Q_n(x)-\lambda_nQ_{n-1}(x),
    \qquad 0\leq n<N,
  \]
  where \( Q_{-1}(x)=0 \), \( Q_0(x)=1 \), and \(\lambda_n\neq0\) for
  \(1\leq n<N\). Let \( \LL \) be the linear functional whose moments
  \(\LL(x^k)=\mu_k\) are given by
  \[
    \sum_{k\geq0}\mu_kz^k
    =\frac{(Q_{N-1}^{(1)})^*(z)}{Q_N^*(z)}.
  \]
  Then \( \{Q_n(x):0\le n\le N \} \) is the monic extended OPS with
  respect to \(\LL\). Moreover, we have
  \[
    \det(\mu_{i+j})_{i,j=0}^{k-1}
    =\prod_{j=1}^{k-1}\lambda_j^{k-j},
    \qquad 1\leq k\leq N.
  \]
\end{prop}

\begin{proof}
  Let \( J \) be the matrix in the proof of
  \Cref{prop:finite-moment-next}. Recall that if
  \( x^n \equiv \sum_{j=0}^{N-1} d_j Q_j \pmod{Q_N} \), then
  \( \LL(x^n)=d_0\). Since \( \LL \) is linear, for any polynomial
  \( p(x) \equiv \sum_{j=0}^{N-1} c_j Q_j \pmod{Q_N} \), we have
  \( \LL(p(x)) = c_0 \). Hence, we have \( \LL(Q_0) = 1 \),
  \( \LL(Q_j) = 0 \) for \( 1\le j\le N-1 \), and
  \( \LL(x^k Q_N(x))=0 \) for \( k\ge0 \). Then, by the same argument
  in the proof of \cite[Chapter~I, Theorem~4.4]{Chihara}, we obtain
  \( \LL(x^i Q_j) = 0 \) and
  \( \LL(x^jQ_j) = \lambda_1 \cdots \lambda_j\ne 0 \) for
  \( 0\le i<j<N \). Therefore, \( \{Q_n(x):0\le n\le N \} \) is the
  monic extended OPS with respect to \(\LL\). The last statement
  follows from \Cref{prop:finite-moment-next}.
\end{proof}

\subsection{Companion OPS}

Let \( f(x) = (x-\alpha_1) \cdots (x-\alpha_N) \) and
\( p_k = \sum_{ j=1 }^{ N } \alpha_j^k \). Since
\[
  \sum_{n\ge0} p_n z^n
  = \frac{(Df)^*(z)}{f^*(z)},
\]
we can rewrite \eqref{eq:12} as
\begin{equation}\label{eq:16}
  \disc(f) = \det\left([z^{i+j}] \frac{(Df)^*(z)}{f^*(z)} \right)_{i,j=0}^{N-1},
\end{equation}
where \([z^r]\) denotes the coefficient of \(z^r\) in the formal power
series expansion at \(z=0\). More generally, the following proposition
expresses the resultant of two polynomials as a Hankel-type
determinant.

\begin{prop}\label{prop:res-hankel}
  Let \( N\ge1 \) and let \(f\) and \(g\) be nonzero polynomials with
  \( \deg f=N \) and \( \deg g=N-1 \). Then
  \begin{equation}\label{eq:res-hankel}
    \res(f,g)
    =(-1)^{\binom N2}\lc(f)^{2N-1}
      \det\left(
        [z^{i+j}]\frac{g^*(z)}{f^*(z)}
      \right)_{i,j=0}^{N-1}.
  \end{equation}
\end{prop}

\begin{proof}
  This follows from combining the equations in \cite[Theorem~3.1(1),
  Corollary~3.1, (3.6), (3.8)]{BikkerUteshev1999}.
\end{proof}

By the same argument in \eqref{eq:17}, if
\( \{Q_n(x)\}= \OPS(\{b_n\},\{\lambda_n\}) \) is an OPS whose moments
satisfy
\[
  \mu_k = \frac{\lc(f)}{\lc(g)} [z^k]\frac{g^*(z)}{f^*(z)} ,\qquad 0\le k\le 2N-2,
\]
then we can evaluate \( \res(f,g) \) using \eqref{eq:res-hankel} and
\eqref{eq:18}. We develop a systematic way to find such an OPS. To do
this we introduce companion OPS.

\begin{defn}\label{def:8}
  Let \(N\geq1\), and let \(f\) and \(g\) be polynomials with
  \(\deg f=N\) and \(\deg g=N-1\). The \emph{\((f,g)\)-companion
    orthogonal polynomial sequence}, or \emph{\((f,g)\)-COPS} for
  short, is the monic extended OPS \(\{Q_n(x):0\leq n\leq N\}\) with
  respect to the linear functional \(\LL\) whose moments
  \(\mu_n=\LL(x^n)\) are given by
  \[
    \sum_{n\ge0} \mu_n z^n =\frac{\lc(f)}{\lc(g)} \frac{g^*(z)}{f^*(z)}.
  \]
\end{defn}

By \Cref{prop:finite-moment-next}, the moments \(\mu_n=\LL(x^n)\) of
the \( (f,g) \)-COPS \(\{Q_n(x):0\leq n\leq N\}\) satisfy
\begin{equation}\label{eq:11}
  \frac{(Q_{N-1}^{(1)})^*(z)}{Q_N^*(z)}
  =\sum_{n\geq0}\mu_nz^n.
\end{equation}

We note that the \( (f,g) \)-COPS does not always exist. If it exists,
then it is unique by \Cref{prop:finite-moment-next}. The following
criterion characterizes COPS by their recurrence coefficients.

\begin{prop}\label{prop:fg-cops-criterion}
  Let \(N\geq1\), let \(f\) and \(g\) be polynomials with
  \(\deg f=N\) and \(\deg g=N-1\), and let
  \(Q=\{Q_n(x):0\leq n\leq N\}\) be a sequence of polynomials.
  Then \(Q\) is an \((f,g)\)-COPS if and only if it satisfies
  \begin{equation}\label{eq:Q-rec}
    Q_{n+1}(x)=(x-b_n)Q_n(x)-\lambda_nQ_{n-1}(x),
    \qquad 0\leq n<N,
  \end{equation}
  where \(Q_{-1}(x)=0\), \(Q_0(x)=1\), and \(\lambda_n\ne0\)
  for \(1\leq n<N\), together with
  \[
    Q_N(x)=\frac{f(x)}{\lc(f)},\qquad
    Q_{N-1}^{(1)}(x)=\frac{g(x)}{\lc(g)}.
  \]
\end{prop}

\begin{proof}
  Suppose that \(Q\) is an \((f,g)\)-COPS. By
  \Cref{prop:finite-moment-next}, it satisfies the required recurrence
  with \(\lambda_n\ne0\). Put \(F=f/\lc(f)\) and \(G=g/\lc(g)\).
  Comparing coefficients in
  \(F^*(z)\sum_{k\geq0}\LL(x^k)z^k=G^*(z)\) gives
  \(\LL(x^kF)=0\) for all \(k\geq0\). Uniqueness in
  \Cref{prop:finite-moment-next} therefore gives \(Q_N=F\), and
  its moment generating function formula then gives
  \(Q_{N-1}^{(1)}=G\).

  Conversely, \Cref{prop:finite-moment-converse} gives a monic
  extended OPS whose moment generating function is
  \((Q_{N-1}^{(1)})^*/Q_N^*=G^*/F^*\), so its moment functional
  is precisely the one in \Cref{def:8}.
\end{proof}

The recurrence coefficients of the \( (f,g) \)-COPS give a product
formula for the resultant of \( f \) and \( g \).

\begin{prop}\label{prop:fg-cops-resultant}
  Let \(N\geq1\), and let \(f\) and \(g\) be polynomials with
  \(\deg f=N\) and \(\deg g=N-1\). Suppose that the \((f,g)\)-COPS
  \(Q=\{Q_n(x):0\leq n\leq N\}\) exists, satisfying the recurrence
  \eqref{eq:Q-rec}. Let \(a=\lc(f)\) and \(b=\lc(g)\). Then
  \[
    \res(f,g)=(-1)^{\binom N2}a^{N-1}b^N
      \prod_{i=1}^{N-1}\lambda_i^{N-i}.
  \]
  In particular, if the \((f,f')\)-COPS exists, then
  \[
    \disc(f)=N^N a^{2N-2}
      \prod_{i=1}^{N-1}\lambda_i^{N-i}.
  \]
\end{prop}
\begin{proof}
  Let \(\{\mu_n:n\geq0\}\) be defined by \eqref{eq:11}. By
  \Cref{prop:res-hankel,prop:finite-moment-converse},
  \[
    \res(Q_N,Q_{N-1}^{(1)})
    =(-1)^{\binom N2}\det(\mu_{i+j})_{i,j=0}^{N-1}
    =(-1)^{\binom N2}
      \prod_{i=1}^{N-1}\lambda_i^{N-i}.
  \]
  Since \( f = a Q_N \) and \( g = b Q_{N-1}^{(1)} \), the definition
  of the resultant gives
  \[
    \res(f,g)=a^{N-1}b^N\res(Q_N,Q_{N-1}^{(1)}),
  \]
  which proves the first formula.

  For \(g=f'\), we have \(b=Na\). Combining the first formula with
  \(\disc(f)=(-1)^{\binom N2}a^{-1}\res(f,f')\) gives the second
  formula.
\end{proof}

\subsection{Construction of COPS}

We now give explicit constructions of \((f_N,Tf_N)\)-COPS and
\((\A_{u,v}f_N,\D_{u,v}f_N)\)-COPS, which provide alternative proofs
of the formulas for the generalized discriminants in terms of their
recurrence coefficients. The following proposition shows that the
mixed three-term relation gives a COPS.

\begin{prop}\label{prop:formula_COPS}
  Suppose that \( \{f_n(x)\} = \OPS(\{b_n\},\{\lambda_n\}) \) and
\( \{\widetilde{f}_n(x)\} =
\OPS(\{\widetilde{b}_n\},\{\widetilde{\lambda}_n\}) \) form a
\( T \)-pair and a linked pair, so that for every \(m\ge1\),
\begin{align}
  \label{eq:14}
  T(f_{m}(x)) &= \eta_m \cdot \widetilde{f}_{m-1}(x),\\
  \label{eq:15}
  f_m(x) &= (x-A_m) \widetilde{f}_{m-1}(x) - B_m \widetilde{f}_{m-2}(x),
\end{align}
for some constants \( A_m \), \( B_m \), and \( \eta_m \). For a fixed
\( N\ge2 \), the \((f_N,Tf_N)\)-COPS \(Q=\{Q_k(x):0\le k\le N\}\) is
given by \( Q_{-1}(x) = 0 \), \(Q_0(x)=1\), and
\[
  Q_{k+1}(x) = (x-b^Q_k) Q_k(x) - \lambda^Q_k Q_{k-1}(x), \qquad 0 \le k<N,
\]
where
\[
b_k^Q =
\begin{cases}
 A_N, & k=0,\\
 \widetilde{b}_{N-1-k}, & 1\le k\le N-1,
\end{cases}
\qquad
\lambda_k^Q =
\begin{cases}
 B_N, & k=1,\\
 \widetilde{\lambda}_{N-k}, & 2\le k\le N-1.
\end{cases}
\]
Here, we assume that \( \lambda_k^Q\ne0 \) for \( 1\le k\le N-1 \).
\end{prop}

\begin{proof}
By the recurrence of \( Q_k(x) \), we have
\begin{equation}\label{eq:10}
  Q_k(x) = 
  \left|
    \begin{matrix}
      x-A_N & -1 & ~ & ~ & ~ & ~\\
      -B_N & x-\widetilde{b}_{N-2} & -1 & ~ & ~ & ~\\
      ~ & -\widetilde{\lambda}_{N-2} & x-\widetilde{b}_{N-3} & \ddots & ~ & ~\\
      ~ & ~ & \ddots & \ddots & \ddots & ~\\
      ~ & ~ & ~ & \ddots & x-\widetilde{b}_{N-k+1} & -1\\
      ~ & ~ & ~ & ~ & -\widetilde{\lambda}_{N-k+1} & x-\widetilde{b}_{N-k}
    \end{matrix}
  \right|.
\end{equation}
  For \(k=N\), deleting the first row and column gives the determinant
  for \(Q_{N-1}^{(1)}(x)\). Reversing its rows and columns and then
  transposing gives the usual tridiagonal determinant for
  \(\widetilde f_{N-1}(x)\). Thus, by \eqref{eq:14},
\[
  Q_{N-1}^{(1)}(x) = \widetilde{f}_{N-1}(x) = \frac{1}{\eta_N} T(f_N(x)) .
\]
Expanding the determinant in \eqref{eq:10} with \( k=N \) along its
first row, reversing the remaining determinants in the same way, and
using \eqref{eq:15}, we obtain
\[
   Q_{N}(x) = (x-A_N) \widetilde{f}_{N-1}(x) - B_N \widetilde{f}_{N-2}(x)= f_N(x).
\]
Since \(f_N\) is monic and \(\lc(Tf_N)=\eta_N\), these identities and
\Cref{prop:fg-cops-criterion} imply that \( Q \) is the \((f_N,Tf_N)\)-COPS.
\end{proof}

\begin{exam}[Laguerre]\label{ex:laguerre-cops}
  Fix \(\alpha>-1\) and \(N\geq2\), and let \(L_n^{(\alpha)}(x)\)
  denote the monic Laguerre polynomials. Their recurrence coefficients
  are
  \[
    b_n = 2n + \alpha + 1, \qquad \lambda_n = n(n+\alpha).
  \]
  The derivative relation
  \[
    D L_n^{(\alpha)}(x)=nL_{n-1}^{(\alpha+1)}(x),
  \]
  and the mixed relation
  \[
    L_n^{(\alpha)}(x)
    =(x-(n+\alpha))L_{n-1}^{(\alpha+1)}(x)
      -(n-1)(n+\alpha)L_{n-2}^{(\alpha+1)}(x)
  \]
  show that \((\{L_n^{(\alpha)}\}, \{L_n^{(\alpha+1)}\})\) is a
  \( D \)-pair and a linked pair, with \(\eta_n=n\), \(A_n=n+\alpha\),
  and \(B_n=(n-1)(n+\alpha)\). By \Cref{prop:formula_COPS}, the
  \((L_N^{(\alpha)},DL_N^{(\alpha)})\)-COPS has recurrence
  coefficients
  \[
    b_k^Q=
    \begin{cases}
      N+\alpha,&k=0,\\
      2N-2k+\alpha,&1\leq k\leq N-1,
    \end{cases}
    \qquad
    \lambda_k^Q=(N-k)(N+\alpha+1-k)
    \quad(1\leq k\leq N-1).
  \]
  Since \( \lc(DL_N^{(\alpha)}) = N \), \Cref{prop:fg-cops-resultant}
  gives
  \[
    \disc(L_N^{(\alpha)})
    =N^N\prod_{k=1}^{N-1}
      \bigl((N-k)(N+\alpha+1-k)\bigr)^{N-k}
    =\prod_{k=1}^N k^k(\alpha+k)^{k-1}.
  \]
\end{exam}

We next construct an \((\A_{u,v}f_N,\D_{u,v}f_N)\)-COPS. The mixed
three-term relation provides the two initial recurrence coefficients.

\begin{prop}\label{prop:two-operator-cops-D-pair}
  Let \(N\geq2\).  Suppose that
  \[
    \{f_n\}=\OPS(\{b_n\},\{\lambda_n\}),
    \qquad
    \{\widetilde f_n\}
      =\OPS(\{\widetilde b_n\},\{\widetilde\lambda_n\})
  \]
  form a \(\D_{u,v}\)-pair, where \(\A_{u,v}x=\rho x+\kappa\). Let
  \( \eta_n \) denote the leading coefficient of \( \D_{u,v}x^n \),
  and let
  \begin{align*}
    A_N
    &=\frac{\eta_{N+1}\widetilde b_{N-1}
      -\eta_N b_N+\kappa\eta_N}{\eta_{N+1}-\rho\eta_N},\\
    B_N
    &=\frac{\eta_{N+1}\widetilde\lambda_{N-1}
      -\eta_{N-1}\lambda_N}{\eta_{N+1}-\rho\eta_N}.
  \end{align*}
  Here, we assume that \(\eta_{N+1}-\rho\eta_N\ne0\) and
  \(B_N\ne0\). Define a polynomial sequence
  \(Q=\{Q_k(x):0\leq k\leq N\}\) by
  \[
    Q_{k+1}(x) = (x-b^Q_k) Q_k(x) - \lambda^Q_k Q_{k-1}(x), \qquad 0\le k<N,
  \]
  where \( Q_{-1}(x) = 0 \), \( Q_0(x) = 1 \), and
  \[
    b_k^Q=
    \begin{cases}
      A_N,&k=0,\\
      \widetilde b_{N-1-k},&1\leq k\leq N-1,
    \end{cases}
    \qquad
    \lambda_k^Q=
    \begin{cases}
      B_N,&k=1,\\
      \widetilde\lambda_{N-k},&2\leq k\leq N-1.
    \end{cases}
  \]
  Then \(Q\) is the \((\A_{u,v}f_N,\D_{u,v}f_N)\)-COPS.
\end{prop}

\begin{proof}
  By \Cref{thm:mixed-3-term},
  \[
    \frac{1}{\eta_{N+1}-\rho\eta_N}\A_{u,v}f_N
    =(x-A_N)\widetilde f_{N-1}
      -B_N\widetilde f_{N-2},
    \qquad
    \frac{1}{\eta_N}\D_{u,v}f_N=\widetilde f_{N-1}.
  \]
  In particular, \( \lc(\A_{u,v}f_N) = \eta_{N+1}-\rho\eta_N \) and
  \( \lc(\D_{u,v}f_N) = \eta_N \). By the same determinant
  construction as in \Cref{prop:formula_COPS}, the displayed
  recurrence coefficients therefore give
  \[
    Q_N=\frac{1}{\eta_{N+1}-\rho\eta_N}\A_{u,v}f_N,
    \qquad
    Q_{N-1}^{(1)}=\widetilde f_{N-1}
      =\frac{1}{\eta_N}\D_{u,v}f_N.
  \]
  Since \(B_N\ne0\) and \(\widetilde\lambda_j\ne0\) for
  \(1\leq j\leq N-2\), all the coefficients \(\lambda_k^Q\) are
  nonzero. Hence \Cref{prop:fg-cops-criterion} applies.
\end{proof}

Under the hypotheses of \Cref{prop:two-operator-cops-D-pair},
\Cref{def:disc2,prop:fg-cops-resultant} give
\[
  \disc(f_N;\A_{u,v},\D_{u,v})
  =\frac{\eta_N^N
    (\eta_{N+1}\widetilde\lambda_{N-1}
      -\eta_{N-1}\lambda_N)^{N-1}}{\eta_{N+1}-\rho\eta_N}
    \prod_{k=1}^{N-2}\widetilde\lambda_k^k,
\]
When \(B_N=0\), the mixed relation shows that
\(\widetilde f_{N-1}\) divides both \(\A_{u,v}f_N\) and
\(\D_{u,v}f_N\), so both sides of this identity vanish. Thus the
identity also holds in this case and recovers \Cref{cor:avg-T-disc}.

\begin{exam}[Askey--Wilson]\label{ex:askey-wilson-cops}
  Fix \(N\geq2\), and let \( p_n(x) = p_n(x;a,b,c,d\mid q)\) be the
  monic Askey--Wilson polynomials, with recurrence coefficients
  \( b_n \) and \( \lambda_n \) given by \eqref{eq:AW-b,la}. Set
  \[
    \widetilde p_n
    =p_n(x;aq^{1/2},bq^{1/2},cq^{1/2},dq^{1/2}\mid q),
  \]
  and let \( \widetilde{b}_n \) and \( \widetilde{\lambda}_n \) denote
  their recurrence coefficients. They are obtained from
  \eqref{eq:AW-b,la} by replacing \( (a,b,c,d) \) with
  \((aq^{1/2},bq^{1/2},cq^{1/2},dq^{1/2})\). The forward shift
  relation is
  \[
    \AW p_n=\eta_n\widetilde p_{n-1},
    \qquad
    \eta_n=\frac{q^{n/2}-q^{-n/2}}
    {q^{1/2}-q^{-1/2}}.
  \]
  Moreover,
  \[
    \A_{\AW}x=\rho x,
    \qquad
    \rho=\frac{q^{1/2}+q^{-1/2}}2,
    \qquad
    \eta_{N+1}-\rho\eta_N
    =\frac{q^{N/2}+q^{-N/2}}2.
  \]
  These formulas yield \( A_N \) and \( B_N \).
  If \(B_N\ne0\), the resulting
  \((\A_{\AW}p_N,\AW p_N)\)-COPS has the recurrence coefficients in
  \Cref{prop:two-operator-cops-D-pair}.
\end{exam}

\subsection{Further remarks on COPS}

We conclude this section with a few remarks on COPS.

For real polynomials, strict interlacing gives a more
transparent criterion and produces a positive finite measure. This is
a first-associated version of Wendroff's theorem; compare
\cite[Chapter~I, Exercise~5.5, p.~30]{Chihara}.

\begin{prop}\label{prop:fg-cops-interlacing}
  Let \(N\ge2\), and let \(f\) and \(g\) be real polynomials of degrees \(N\) and
  \(N-1\), respectively. There exists a unique \((f,g)\)-COPS
    \(\{Q_n(x):0\leq n\leq N\}\) satisfying
  \begin{equation}\label{eq:21}
    Q_{n+1}(x) = (x-b_n) Q_n(x) - \lambda_n Q_{n-1}(x), \qquad 0\le n<N
  \end{equation}
  where \( Q_{-1}(x)=0 \) and \( Q_0(x)=1 \), with real \(b_n\) and
  positive \(\lambda_n\) if and only if \(f\) and \(g\) have real,
  simple, strictly interlacing zeros; that is, the zeros
  \( \alpha_1,\dots,\alpha_N \) of \( f \) and the zeros
  \( \beta_1,\dots,\beta_{N-1} \) of \( g \) satisfy
  \[
    \alpha_1 < \beta_1 <\alpha_2 < \beta_2<\cdots < \alpha_{N-1}< \beta_{N-1}<\alpha_N.
  \]

  In this case, the linear functional \( \LL \) of this COPS is
  represented by
  \begin{equation}\label{eq:fg-cops-measure}
    \LL(x^n)=\sum_{j=1}^N w_j\alpha_j^n,
    \qquad
    w_j=\frac{\lc(f)}{\lc(g)}\frac{g(\alpha_j)}{f'(\alpha_j)},
    \qquad n\ge0.
  \end{equation}
\end{prop}

\begin{proof}
  Suppose first that the zeros of \( f \) and \( g \) are real,
  simple, and strictly interlace. Partial fractions give
  \begin{equation}\label{eq:fg-partial-fractions}
    \frac{\lc(f)}{\lc(g)}\frac{g(z)}{f(z)}
    =\sum_{j=1}^N\frac{w_j}{z-\alpha_j},
    \qquad
    w_j=\frac{\lc(f)}{\lc(g)}\frac{g(\alpha_j)}{f'(\alpha_j)}.
  \end{equation}
  Strict interlacing gives \(w_j>0\), and comparison at infinity
  gives \(\sum_{j=1}^N w_j=1\). Thus
  \eqref{eq:fg-cops-measure} defines a probability measure on \( N \)
  distinct points, whose moment generating function is
  \(\frac{\lc(f)}{\lc(g)}\frac{g^*(t)}{f^*(t)}\) by
  \eqref{eq:fg-partial-fractions}.
  Its Hankel determinants satisfy \(\Delta_n>0\) for \(1\le n\le N\)
  (see \cite[Chapter~I, (3.7), p.~16, and Exercise~5.1, p.~30]{Chihara}).
  Hence \Cref{prop:finite-moment-next} gives the unique
  \((f,g)\)-COPS, with real \(b_j\) and
  \(\lambda_j=\Delta_{j+1}\Delta_{j-1}/\Delta_j^2>0\).

  Conversely, suppose that the \( (f,g) \)-COPS has real \(b_n\) and
  positive \(\lambda_n\). By \Cref{prop:fg-cops-criterion}, we have
  \(f=\lc(f)Q_N\) and \(g=\lc(g)Q_{N-1}^{(1)}\). Extend the recurrence
  \eqref{eq:21} by choosing \(b_n\in\mathbb R\) and \(\lambda_n>0\)
  for \(n\ge N\). By Favard's theorem, the resulting infinite sequence
  \( \{Q_n:n\ge0\} \) is an OPS with respect to a new
  positive-definite linear functional, with the same \(Q_N\) and
  \(Q_{N-1}^{(1)}\). The interlacing theorem for first-associated
  polynomials \cite[Chapter~III, Theorem~4.1, p.~86]{Chihara} shows
  that their zeros are real, simple, and strictly interlace.
\end{proof}

Let \( \{f_n\}= \OPS(\{b_n\},\{\lambda_n\}) \) with real \( b_n \) and
positive \( \lambda_n \). Fix \( N\ge2 \) and let
\( \alpha_1, \dots, \alpha_N \) be the zeros of \( f_N(x) \). The
zeros of \( f_N(x) \) interlace with those of \( f_N'(x) \) by Rolle's
theorem. They also interlace with those of \( f_{N-1}(x) \) (see
\cite[Ch.~I, Theorem~5.3]{Chihara}) and \( f^{(1)}_{N-1}(x) \) (see
\cite[(3.2)]{Chihara1957}). Thus, by \Cref{prop:fg-cops-interlacing},
there exist a unique \( (f_N,g) \)-COPS for
\( g\in \{f_N', f^{(1)}_{N-1}, f_{N-1}\} \) with linear functional
\( \LL \) given by
\[
    \LL(x^n)=\sum_{j=1}^N w_j\alpha_j^n,
    \qquad
    w_j=\frac{1}{\lc(g)}\frac{g(\alpha_j)}{f'_N(\alpha_j)},
    \qquad n\ge0.
\]
We investigate each COPS and the weight \( w_j \) of the measure at
\( \alpha_j \).

\begin{exam}[\((f_N,f'_N)\)-COPS]
  When \(g=f'_N\), we have \(w_j=1/N\). Hence the linear functional of
  the \((f_N,f'_N)\)-COPS is given by
\[
  \LL(x^n)=\frac1N\sum_{j=1}^N \alpha_j^n, \qquad n\ge0.
\]
The representing measure is the normalized zero-counting measure of
\(f_N\). See \cite[Theorem~2.4]{Simon2009} for its relation to weak
limits of Christoffel--Darboux kernel measures.
\end{exam}

\begin{exam}[\((f_N,f_{N-1}^{(1)})\)-COPS]
The first associated polynomials satisfy the Casoratian identity
(see \cite[Ch.~III, (4.4), p.~86]{Chihara})
\[
  f_{N-1}(x)f_{N-1}^{(1)}(x)
  -f_N(x)f_{N-2}^{(1)}(x)
  =\lambda_1\cdots\lambda_{N-1}.
\]
Thus, at a zero \(\alpha_j\) of \(f_N\), we have
\( f_{N-1}^{(1)}(\alpha_j) ={\lambda_1\cdots\lambda_{N-1}}/{f_{N-1}(\alpha_j)} \), which implies
\begin{equation}\label{eq:20}
  w_j
  =\frac{f_{N-1}^{(1)}(\alpha_j)}{f_N'(\alpha_j)}
  =\frac{\lambda_1\cdots\lambda_{N-1}} {f_N'(\alpha_j)f_{N-1}(\alpha_j)}.
\end{equation}
These are precisely the weights for the Gaussian quadrature formula
for the linear functional \(\LL_f\) of the OPS \( \{f_n\} \),
normalized by \(\LL_f(1)=1\) (see \cite[Ch.~I,
Exercise~6.2(a), p.~34]{Chihara}). Therefore, the
\((f_N,f_{N-1}^{(1)})\)-COPS is the original OPS
\(\{f_n:0\leq n\leq N\}\) truncated at degree \(N\).
\end{exam}

\begin{exam}[\((f_N,f_{N-1})\)-COPS]
  The finite measure of the \((f_N,f_{N-1})\)-COPS has weights
\[
  w_j
  =\frac{f_{N-1}(\alpha_j)}{f_N'(\alpha_j)}
  =\omega_j\frac{f_{N-1}(\alpha_j)^2}{\lambda_1\cdots\lambda_{N-1}},
\]
where \( \omega_j \) is the right-hand side of \eqref{eq:20}. Thus this
measure is obtained from the Gaussian quadrature measure by multiplying
each weight \(\omega_j\) at \(\alpha_j\) by
\(f_{N-1}(\alpha_j)^2/\lambda_1\cdots\lambda_{N-1}\).
These factors are positive, since strict interlacing gives
\(f_{N-1}(\alpha_j)\ne0\) and all \(\lambda_k\) are positive.
\end{exam}

\section{Final remarks}
\label{sec:final-remarks}

In this paper, we introduced the generalized discriminant
\(\disc(f;\A_{u,v},\D_{u,v})\) of a polynomial \(f\), defined with
respect to the averaging operator \(\A_{u,v}\) and the divided
difference operator \(\D_{u,v}\). It extends the ordinary discriminant
while retaining analogous properties. For example,
\Cref{pro:disc-prod2} gives a product formula in terms of the zeros of
\(f\), and our results yield product formulas for the corresponding
generalized discriminants of every OPS in the Askey and \(q\)-Askey
schemes. It would be interesting to determine whether other features
of ordinary discriminants, such as their electrostatic interpretation,
have analogues in this generalized setting.

In the original electrostatics problem, see
\cite[Section~3.5]{Ismail}, $n$ unit charges are placed inside
$[-1,1]$ at locations $-1<x_1<x_2<\cdots<x_n<1$. Two charges,
depending on $\alpha$ and $\beta$, are placed at $-1$ and $1$. A
potential involving the discriminant is applied. The equilibrium
locations of the particles are the zeros of the \(n\)th Jacobi
polynomial \(P_n^{(\alpha,\beta)}(x)\). One
may ask if there is a similar electrostatics problem, with a resultant
replacing the discriminant, and the equilibrium locations being the
interlaced zeros of consecutive Jacobi polynomials. The answer is yes,
and involves two types of particles, located at
$-1<x_1<y_1<x_2<y_2<\cdots<y_{n-1}<x_n<1$ which have no
self-interaction. We do not make a formal statement of this result in
this paper.

The suitably normalized Askey--Wilson polynomials
\(c_n p_n(x;a,bq^{-n},c,d\mid q)\) are type \(R_{II}\) polynomials in
both \(a\) and \(b\). Thus, by the type \(R_{II}\) version of Schur's
lemma (\Cref{SchurRII}), their resultants with respect to \(a\) and
\(b\) factor. These factorizations may be related to moment duality
\cite{Ismail2025}.

\appendix

\section{Generalized Discriminant Product Formulas}
\label{sec:app}

Each orthogonal polynomial sequence
\( \{p_n(x)\} = \OPS(\{b_n\},\{\lambda_n\}) \) in the Askey and
\( q \)-Askey schemes has a forward shift operator \( T \) such that
the monic polynomials \(\widetilde p_n=Tp_{n+1}/\eta_{n+1}\) form an OPS
\( \{\widetilde p_n(x)\} = \OPS(\{\widetilde{b}_n\},\{\widetilde{\lambda}_n\}) \),
where \(\eta_n\) is the leading coefficient of \(Tx^n\).
We refer the reader to \cite{KLS} for the precise
description of the recurrence coefficients
\( b_n, \lambda_n, \widetilde{b}_n, \widetilde{\lambda}_n\) and the
forward shift operators \( T \).

In this appendix, we provide a table for the monic orthogonal
polynomials \( p_n(x) \) in the Askey and \( q \)-Askey schemes
together with their forward shift operators \( T \) and the quantities
\(\eta_{n+1}\widetilde{\lambda}_{n-1}-\eta_{n-1}\lambda_n\). For
completeness and to make this appendix self-contained, we also record
\( \widetilde{\lambda}_n \). For each family, the upper and lower rows
display \(\eta_{n+1}\widetilde{\lambda}_{n-1}-\eta_{n-1}\lambda_n\)
and \(\widetilde{\lambda}_n\), respectively; most of which have simple
product formulas. Using \Cref{cor:avg-T-disc,cor:linear-disc-ratio}
and \Cref{tab:askey-duv-constants}, one can evaluate
\( \disc(p_n(x);\A_T,T) \). In the affine case
\(u(x)=\alpha x+\beta\), \(v(x)=x+\gamma\), one can also evaluate
\(\disc(p_n(x);S)\), where \(Sf=(Tf)\circ v^{-1}\); in particular,
\(S=T\) when \(v(x)=x\). For \(D_{q^{-1}}\), use the \(D_q\) row of
\Cref{tab:askey-duv-constants} with \(q\) replaced by \(q^{-1}\). Here
\(s=a+b+c+d\) for Wilson and continuous Hahn polynomials. In the first
column, C. and D. abbreviate continuous and discrete, respectively.

\begingroup
\footnotesize
\setlength{\tabcolsep}{1pt}
\renewcommand{\arraystretch}{1.35}
\newcommand{\oplink}[2]{#2}
\newcommand{\ShiftFormulaPair}[2]{%
  \begin{tabular}[t]{@{}p{\linewidth}@{}}
    \({\scriptstyle #1}\)\\[3pt]
    \noalign{%
      \hbox to\linewidth{%
        \leaders\hbox{\rule{\arrayrulewidth}{\arrayrulewidth}\hskip1.5pt}\hfill}%
      \vskip3pt}
    \({\scriptstyle #2}\)
  \end{tabular}}
\newcommand{\ShiftSummaryRow}[4]{#1 & \(#2\) &
  \ShiftFormulaPair{#3}{#4}\\[5pt]
  \hline\noalign{\vskip5pt}}
\begin{longtable}{@{}>{\raggedright\arraybackslash}p{0.18\linewidth}
  >{\centering\arraybackslash}p{0.095\linewidth}
  >{\raggedright\arraybackslash}p{0.72\linewidth}@{}}
\label{tab:askey-shifted-lambda}\\
\textbf{\(p_n(x)\)} & \textbf{\(T\)} &
\ShiftFormulaPair{\eta_{n+1}\widetilde{\lambda}_{n-1}-\eta_{n-1}\lambda_n}{\widetilde{\lambda}_n}\\[5pt]
\hline\noalign{\vskip5pt}
\endfirsthead
\endhead
\ShiftSummaryRow{\oplink{ops:wilson}{Wilson}}{\mathcal{W}}
  {\tfrac{ (n-1)(n+a+b-1)(n+a+c-1) (n+a+d-1)(n+b+c-1)(n+b+d-1) (n+c+d-1) }{(2n+s-3)(2n+s-2)^2}}
  {\tfrac{ n(n+s)(n+a+b)(n+a+c) (n+a+d)(n+b+c)(n+b+d) (n+c+d) }{(2n+s-1)(2n+s)^2(2n+s+1)}}
\ShiftSummaryRow{\oplink{ops:racah}{Racah}}{\mathcal{R}}
  {\tfrac{ (n-1)(n+\alpha)(n+\beta) (n+\gamma)(n+\beta+\delta)(n+\alpha-\delta) (n+\alpha+\beta-\gamma) }{(2n+\alpha+\beta-1)(2n+\alpha+\beta)^2}}
  {\tfrac{ n(n+\alpha+\beta+2)(n+\alpha+1)(n+\beta+\delta+1) (n+\gamma+1)(n+\alpha+\beta-\gamma+1)(n+\alpha-\delta+1) (n+\beta+1) }{(2n+\alpha+\beta+1)(2n+\alpha+\beta+2)^2(2n+\alpha+\beta+3)}}
\ShiftSummaryRow{\oplink{ops:continuous-dual-hahn}{C. dual Hahn}}{\mathcal{W}}
  {(n-1)(n+a+b-1)(n+a+c-1) (n+b+c-1)}
  {n(n+a+b)(n+a+c)(n+b+c)}
\ShiftSummaryRow{\oplink{ops:continuous-hahn}{C. Hahn}}{\delta}
  {\tfrac{ (n-1)(n+a+c-1)(n+a+d-1) (n+b+c-1)(n+b+d-1) }{(2n+s-3)(2n+s-2)^2}}
  {\tfrac{ n(n+s)(n+a+c)(n+a+d) (n+b+c)(n+b+d) }{(2n+s-1)(2n+s)^2(2n+s+1)}}
\ShiftSummaryRow{\oplink{ops:hahn}{Hahn}}{\Deltaop}
  {\tfrac{ (n-1)(n+\alpha)(n+\beta) (N-n+1)(n+\alpha+\beta+N+1) }{(2n+\alpha+\beta-1)(2n+\alpha+\beta)^2}}
  {\tfrac{ n(n+\alpha+\beta+2)(n+\alpha+1)(N-n) (n+\alpha+\beta+N+2)(n+\beta+1) }{(2n+\alpha+\beta+1)(2n+\alpha+\beta+2)^2(2n+\alpha+\beta+3)}}
\ShiftSummaryRow{\oplink{ops:dual-hahn}{Dual Hahn}}{\mathcal{R}}
  {(n-1)(n+\gamma)(n-N-1) (n-\delta-N-1)}
  {n(n+\gamma+1)(n-N)(n-\delta-N)}
\ShiftSummaryRow{\oplink{ops:meixner-pollaczek}{Meixner--Pollaczek}}{\delta}
  {\tfrac{(n-1)(n+2\lambda-1)}{4\sin^2\phi}}
  {\tfrac{n(n+2\lambda)}{4\sin^2\phi}}
\ShiftSummaryRow{\oplink{ops:jacobi}{Jacobi}}{D}
  {\tfrac{4(n-1)(n+\alpha)(n+\beta)}{(2n+\alpha+\beta-1)(2n+\alpha+\beta)^2}}
  {\tfrac{4n(n+\alpha+1)(n+\beta+1)(n+\alpha+\beta+2)}{(2n+\alpha+\beta+1)(2n+\alpha+\beta+2)^2(2n+\alpha+\beta+3)}}
\ShiftSummaryRow{\oplink{ops:meixner}{Meixner}}{\Deltaop}
  {(n-1)(n+\beta-1)c(1-c)^{-2}}
  {n(n+\beta)c(1-c)^{-2}}
\ShiftSummaryRow{\oplink{ops:krawtchouk}{Krawtchouk}}{\Deltaop}
  {(n-1)p(1-p) (N+1-n)}
  {np(1-p)(N-n)}
\ShiftSummaryRow{\oplink{ops:laguerre}{Laguerre}}{D}
  {(n-1)(n+\alpha)}
  {n(n+\alpha+1)}
\ShiftSummaryRow{\oplink{ops:bessel}{Bessel}}{D}
  {-\tfrac{4(n-1)}{(2n+a-1)(2n+a)^2}}
  {-\tfrac{4n(n+a+2)}{(2n+a+1)(2n+a+2)^2(2n+a+3)}}
\ShiftSummaryRow{\oplink{ops:charlier}{Charlier}}{\Deltaop}
  {(n-1)a}
  {na}
\ShiftSummaryRow{\oplink{ops:hermite}{Hermite}}{D}
  {\tfrac{n-1}{2}}
  {\tfrac{n}{2}}
\ShiftSummaryRow{\oplink{ops:q-askey-wilson}{Askey--Wilson}}{\AW}
  {\tfrac{ q^{-n/2}(1-q^{n-1})(1-abq^{n-1})(1-acq^{n-1}) (1-adq^{n-1})(1-bcq^{n-1})(1-bdq^{n-1}) (1-cdq^{n-1}) }{4(1-abcdq^{2n-3})(1-abcdq^{2n-2})^2}}
  {\tfrac{ (1-q^n)(1-abcdq^n)(1-abq^n) (1-acq^n)(1-adq^n)(1-bcq^n) (1-bdq^n)(1-cdq^n) }{4(1-abcdq^{2n-1})(1-abcdq^{2n})^2(1-abcdq^{2n+1})}}
\ShiftSummaryRow{\oplink{ops:q-racah}{\(q\)-Racah}}{\Rq}
  {\tfrac{ q(q^n-q)(1-\alpha q^n)(1-\beta q^n) (1-\gamma q^n)(1-\beta\delta q^n)(\gamma-\alpha\beta q^n) (\delta-\alpha q^n) }{q^{n/2}(\alpha\beta q^{2n}-q)(\alpha\beta q^{2n}-1)^2}}
  {\tfrac{ q(1-q^n)(1-\alpha q^{n+1})(1-\beta q^{n+1}) (1-\gamma q^{n+1})(1-\beta\delta q^{n+1})(1-\alpha\beta q^{n+2}) (\gamma-\alpha\beta q^{n+1})(\delta-\alpha q^{n+1}) }{(1-\alpha\beta q^{2n+1})(1-\alpha\beta q^{2n+2})^2(1-\alpha\beta q^{2n+3})}}
\ShiftSummaryRow{\oplink{ops:q-continuous-dual-hahn}{C. dual \(q\)-Hahn}}{\AW}
  {\tfrac{ 1 }{4}q^{-n/2}(1-q^{n-1})(1-abq^{n-1})(1-acq^{n-1}) (1-bcq^{n-1})}
  {\tfrac{ 1 }{4}(1-q^n)(1-abq^n)(1-acq^n) (1-bcq^n)}
\ShiftSummaryRow{\oplink{ops:q-continuous-hahn}{C. \(q\)-Hahn}}{\AW}
  {\tfrac{ q^{-n/2}(1-q^{n-1})(1-abe^{2\ii\phi}q^{n-1})(1-acq^{n-1}) (1-adq^{n-1})(1-bcq^{n-1})(1-bdq^{n-1}) (1-cde^{-2\ii\phi}q^{n-1}) }{4(1-abcdq^{2n-3})(1-abcdq^{2n-2})^2}}
  {\tfrac{ (1-q^n)(1-abcdq^n)(1-abe^{2\ii\phi}q^n) (1-acq^n)(1-adq^n)(1-bcq^n) (1-bdq^n)(1-cde^{-2\ii\phi}q^n) }{4(1-abcdq^{2n-1})(1-abcdq^{2n})^2(1-abcdq^{2n+1})}}
\ShiftSummaryRow{\oplink{ops:q-big-jacobi}{Big \(q\)-Jacobi}}{D_q}
  {-\tfrac{ acq^n(1-q^{n-1})(1-aq^n)(1-bq^n) (1-cq^n)(1-abc^{-1}q^n) }{(1-abq^{2n-1})(1-abq^{2n})^2}}
  {-\tfrac{ acq^{n+1}(1-q^n)(1-aq^{n+1})(1-bq^{n+1}) (1-cq^{n+1})(1-abq^{n+2})(1-abc^{-1}q^{n+1}) }{(1-abq^{2n+1})(1-abq^{2n+2})^2(1-abq^{2n+3})}}
\ShiftSummaryRow{\oplink{ops:q-hahn}{\(q\)-Hahn}}{D_{q^{-1}}}
  {-\tfrac{ \alpha q^{1-N}(1-q^{n-1})(1-q^{n-N-1})(1-\alpha q^n) (1-\beta q^n)(1-\alpha\beta q^{n+N+1}) }{(1-\alpha\beta q^{2n-1})(1-\alpha\beta q^{2n})^2}}
  {-\tfrac{ \alpha q^{n-N+2}(1-q^n)(1-q^{n-N})(1-\alpha q^{n+1}) (1-\beta q^{n+1})(1-\alpha\beta q^{n+2})(1-\alpha\beta q^{n+N+2}) }{(1-\alpha\beta q^{2n+1})(1-\alpha\beta q^{2n+2})^2(1-\alpha\beta q^{2n+3})}}
\ShiftSummaryRow{\oplink{ops:q-dual-hahn}{Dual \(q\)-Hahn}}{\Rq}
  {\gamma q^{1-n/2}(1-q^{n-1})(1-q^{n-N-1})(1-\gamma q^n)(\delta-q^{n-N-1})}
  {\gamma q(1-q^n)(1-q^{n-N})(1-\gamma q^{n+1})(\delta-q^{n-N})}
\ShiftSummaryRow{\oplink{ops:q-al-salam-chihara}{Al-Salam--Chihara}}{\AW}
  {\tfrac{1}{4}q^{-n/2}(1-q^{n-1})(1-abq^{n-1})}
  {\tfrac{1}{4}(1-q^n)(1-abq^n)}
\ShiftSummaryRow{\oplink{ops:q-meixner-pollaczek}{\(q\)-Meixner--Pollaczek}}{\AW}
  {\tfrac{1}{4}q^{-n/2}(1-q^{n-1})(1-a^2q^{n-1})}
  {\tfrac{1}{4}(1-q^n)(1-a^2q^n)}
\ShiftSummaryRow{\oplink{ops:q-continuous-jacobi}{C. \(q\)-Jacobi}}{\AW}
  {\tfrac{ (1-q^{n-1})(1-q^{n+\alpha})(1-q^{n+\beta}) (1+q^{n+(\alpha+\beta+1)/2}) }{4q^{n/2}(1-q^{n+(\alpha+\beta-1)/2}) (1-q^{n+(\alpha+\beta)/2})^2}}
  {\tfrac{ (1-q^n)(1-q^{n+\alpha+1})(1-q^{n+\beta+1}) (1-q^{n+\alpha+\beta+2}) }{4(1-q^{n+(\alpha+\beta+1)/2})(1-q^{n+(\alpha+\beta+2)/2})^2(1-q^{n+(\alpha+\beta+3)/2})}}
\ShiftSummaryRow{\oplink{ops:q-big-laguerre}{Big \(q\)-Laguerre}}{D_q}
  {-abq^n(1-q^{n-1})(1-aq^n)(1-bq^n)}
  {-abq^{n+1}(1-q^n)(1-aq^{n+1})(1-bq^{n+1})}
\ShiftSummaryRow{\oplink{ops:q-little-jacobi}{Little \(q\)-Jacobi}}{D_q}
  {\tfrac{aq^{2n-2}(1-q^{n-1})(1-aq^n)(1-bq^n)}{(1-abq^{2n-1})(1-abq^{2n})^2}}
  {\tfrac{ aq^{2n}(1-q^n)(1-aq^{n+1})(1-bq^{n+1}) (1-abq^{n+2}) }{(1-abq^{2n+1})(1-abq^{2n+2})^2(1-abq^{2n+3})}}
\ShiftSummaryRow{\oplink{ops:q-meixner}{\(q\)-Meixner}}{D_{q^{-1}}}
  {cq^{2-4n}(q-q^n)(1-bq^n)(c+q^n)}
  {cq^{-4n-1}(1-q^n)(1-bq^{n+1})(c+q^{n+1})}
\ShiftSummaryRow{\oplink{ops:q-quantum-krawtchouk}{Quantum \(q\)-Krawtchouk}}{D_{q^{-1}}}
  {p^{-2}q^{2-4n}(q-q^n)(1-pq^n)(1-q^{n-N-1})}
  {p^{-2}q^{-4n-1}(1-q^n)(1-pq^{n+1})(1-q^{n-N})}
\ShiftSummaryRow{\oplink{ops:q-krawtchouk}{\(q\)-Krawtchouk}}{D_{q^{-1}}}
  {-\tfrac{pq^{n-N}(1-q^{n-1})(1-q^{n-N-1})(1+pq^{n+N})}{(1+pq^{2n-2})(1+pq^{2n-1})^2}}
  {-\tfrac{pq^{2n-N+2}(1-q^n)(1-q^{n-N})(1+pq^{n+1})(1+pq^{n+N+1})}{(1+pq^{2n})(1+pq^{2n+1})^2(1+pq^{2n+2})}}
\ShiftSummaryRow{\oplink{ops:q-affine-krawtchouk}{Affine \(q\)-Krawtchouk}}{D_{q^{-1}}}
  {-pq^{1-N}(1-q^{n-1})(1-pq^n)(1-q^{n-N-1})}
  {-pq^{n-N+2}(1-q^n)(1-pq^{n+1})(1-q^{n-N})}
\ShiftSummaryRow{\oplink{ops:q-dual-krawtchouk}{Dual \(q\)-Krawtchouk}}{\Rq}
  {cq^{-N-n/2}(1-q^{n-1})(1-q^{n-N-1})}
  {cq^{-N}(1-q^n)(1-q^{n-N})}
\ShiftSummaryRow{\oplink{ops:q-continuous-big-hermite}{C. big \(q\)-Hermite}}{\AW}
  {\tfrac{1}{4}q^{-n/2}(1-q^{n-1})}
  {\tfrac{1}{4}(1-q^n)}
\ShiftSummaryRow{\oplink{ops:q-continuous-laguerre}{C. \(q\)-Laguerre}}{\AW}
  {\tfrac{1}{4}q^{-n/2}(1-q^{n-1})(1-q^{n+\alpha})}
  {\tfrac{1}{4}(1-q^n)(1-q^{n+\alpha+1})}
\ShiftSummaryRow{\oplink{ops:q-little-laguerre}{Little \(q\)-Laguerre}}{D_q}
  {aq^{2n-2}(1-q^{n-1})(1-aq^n)}
  {aq^{2n}(1-q^n)(1-aq^{n+1})}
\ShiftSummaryRow{\oplink{ops:q-laguerre}{\(q\)-Laguerre}}{D_q}
  {q^{-3n-2\alpha+1}(1-q^{n-1})(1-q^{n+\alpha})}
  {q^{-4n-2\alpha-3}(1-q^n)(1-q^{n+\alpha+1})}
\ShiftSummaryRow{\oplink{ops:q-bessel}{\(q\)-Bessel}}{D_q}
  {\tfrac{aq^{3n-3}(1-q^{n-1})}{(1+aq^{2n-2})(1+aq^{2n-1})^2}}
  {\tfrac{aq^{3n}(1-q^n)(1+aq^{n+1})}{(1+aq^{2n})(1+aq^{2n+1})^2(1+aq^{2n+2})}}
\ShiftSummaryRow{\oplink{ops:q-charlier}{\(q\)-Charlier}}{D_{q^{-1}}}
  {aq^{2-4n}(q-q^n)(a+q^n)}
  {aq^{-4n-1}(1-q^n)(a+q^{n+1})}
\ShiftSummaryRow{\oplink{ops:q-al-salam-carlitz-i}{Al-Salam--Carlitz~I}}{D_q}
  {-aq^{n-2}(1-q^{n-1})}
  {-aq^{n-1}(1-q^n)}
\ShiftSummaryRow{\oplink{ops:q-al-salam-carlitz-ii}{Al-Salam--Carlitz~II}}{D_q}
  {aq^{-n+1}(1-q^{n-1})}
  {aq^{-2n-1}(1-q^n)}
\ShiftSummaryRow{\oplink{ops:q-continuous-hermite}{C. \(q\)-Hermite}}{\AW}
  {\tfrac{1}{4}q^{-n/2}(1-q^{n-1})}
  {\tfrac{1}{4}(1-q^n)}
\ShiftSummaryRow{\oplink{ops:q-stieltjes-wigert}{Stieltjes--Wigert}}{D_q}
  {q^{-3n+1}(1-q^{n-1})}
  {q^{-4n-3}(1-q^n)}
\ShiftSummaryRow{\oplink{ops:q-discrete-hermite-i}{D. \(q\)-Hermite~I}}{D_q}
  {q^{n-2}(1-q^{n-1})}
  {q^{n-1}(1-q^n)}
\ShiftSummaryRow{\oplink{ops:q-discrete-hermite-ii}{D. \(q\)-Hermite~II}}{D_q}
  {q^{-n+1}(1-q^{n-1})}
  {q^{-2n-1}(1-q^n)}
\end{longtable}
\endgroup

\bibliographystyle{abbrv}

\end{document}